\documentclass[preprint,9pt]{elsarticle}

 \usepackage{graphics}

\usepackage{amssymb}

\newtheorem{theorem}{Theorem}
\newtheorem{lemma}[theorem]{Lemma}

\newtheorem{remark}[theorem]{Remark}

\journal{Computational Physics}

\begin{document}

\begin{frontmatter}



\title{Efficient determination of the critical parameters and the statistical quantities for Klein-Gordon and sine-Gordon equations with a singular potential using generalized polynomial chaos methods}

\author{ Debananda Chakraborty and Jae-Hun Jung*}
\ead{dc58@buffalo.edu, jaehun@buffalo.edu}
\cortext[cor1]{Corresponding author}


\address{Department of Mathematics, State University of New York at Buffalo, Buffalo, NY 14260-2900, USA.}

\begin{abstract}
We consider the Klein-Gordon and sine-Gordon type equations with a point-like potential, which describes the wave phenomenon in disordered media with a defect. The singular potential term yields a critical phenomenon--that is, the solution behavior around the critical parameter value bifurcates into two extreme cases. Finding such critical parameter values and the associated statistical quantities demands a large number of individual simulations with different parameter values. Pinpointing the critical value with arbitrary accuracy is even more challenging. In this work, we adopt the generalized polynomial chaos (gPC) method to determine the critical values and the mean solutions around such values. 

First, we consider the critical value associated with the strength of the singular potential for the Klein-Gordon equation. We show the existence of a critical behavior with certain boundary conditions. Then we expand the solution in the random variable associated with the parameter. The obtained partial differential equations are solved using the Chebyshev collocation method. Due to the existence of the singularity, the Gibbs phenomenon appears in the solution, yielding a slow convergence of the numerically computed critical value. To deal with the singularity, we adopt the consistent spectral collocation method. The gPC method, along with the consistent Chebyshev method, determines the critical value and the mean solution highly efficiently. 

We then consider the sine-Gordon equation, for which the critical value is associated with the initial velocity of the kink soliton solution. The critical behavior in this case is that the solution \textit{passes} through (particle-pass), is \textit{trapped} by (particle-capture), or is \textit{reflected} by (particle-reflection) the singular potential if the initial velocity of the soliton solution is greater than, equal to, or less than the critical value, respectively. Due to the nonlinearity of the equation, we use the gPC mean value rather than reconstructing the solution to find the critical parameter. Numerical results show that the critical value can be determined efficiently and accurately by using the proposed method. The results are also compared with the results using the Monte-Carlo method.

\end{abstract}


\end{frontmatter}
\section{Introduction}
Critical phenomena of partial differential equations (PDEs) associated with the critical parameter can be found in various physical applications. For example, in \cite{Choptuik} the gravitational collapse of the scalar field was considered. The critical phenomenon of the gravitational collapse was discovered; there exists a critical initial mass such that if the initial mass is less than the critical mass, the scalar field eventually disperses. Conversely, it forms a black hole if the value exceeds the critical mass. If the initial parameter is arbitrarily close to the critical value,  scale-invariant behavior of the solution occurs \cite{Choptuik}. For PDEs from nonlinear optics, a similar phenomenon has been found and studied, e.g. the sine-Gordon equation or nonlinear Schr\"{o}dinger equation with a point-like singular potential term \cite{Fei,GH2004,GHW2002,PZ2007}. For the sine-Gordon equation, there exists a kink soliton solution; this yields a critical behavior once a self-interacting singular potential is placed. This phenomenon is known as the wave phenomenon in disordered media with a defect. If the initial velocity of the kink soliton solution exceeds the critical velocity, the kink soliton eventually passes through the potential, which is known as the {\it particle-pass}. If the initial velocity is less than the critical velocity, the soliton is either captured by or reflected by the potential, known as the {\it particle-capture} or {\it particle-reflection}, respectively. Since the sine-Gordon equation is not integrable when the singular potential term is present, no exact solution is available; no exact value of the critical parameter is available either. There are many unknown properties of the solution behaviors in this case. For more detailed study, we refer the readers to \cite{Fei,GH2004,GH2007,GHW2002,GSW2002}. A similar phenomenon has been also found with nonlinear Schr\"{o}dinger equations \cite{MMC2003,TSB2003}. For example, in \cite{CM1995,MMC2003} the singular potential term perturbs the soliton propagation. Similar phenomena of particle-pass, particle-capture, and particle-reflection were observed for some critical parameters. 

In this paper, we consider the critical phenomenon that is induced by the point-like singular source term in the PDEs, such as the sine-Gordon equation mentioned above. In particular, we consider the Klein-Gordon and sine-Gordon equations to be perturbed by the singular potential term. For both cases, the solution behaviors are highly sensitive to the critical value of the parameter. Determination of the critical value and the mean solution is an important task because the wave or soliton interactions with the local defect should be well understood for better optical communications. In either particle-capture or particle-reflection, the optical signal is not transmitted with the associated parametric value if it is less than the critical value. The signal should be transmitted with some values of the parameter that exceed the critical value for the communication without blocking or trapping \cite{EB2010,MV2006}. Finding the critical value is doable but challenging. 

There are two reasons. First, the singular potential is given by the Dirac $\delta$-function, so the given PDEs are not well-defined. For the sine-Gordon equation, no exact solution is known with the presence of the singular potential term, which means that no analytic form of the critical value is available--only the numerical or real laboratory experimental approach can render the critical value measurable. The numerical determination of the critical parameter value relies on expensive Monte-Carlo(MC) simulations based on the bi-section method. Moreover, it is highly expensive to determine the critical value up to an arbitrary accuracy (for the sine-Gordon case, it is not even known whether the critical value is rational or irrational). Furthermore, the numerically obtained critical value depends on the numerical method, considered with given numerical settings such as the truncation number, domain size, final time  etc. \cite{Wang}. Due to the singularity in the potential term, the Gibbs phenomenon also appears in the finite truncated solution, resulting in a slow convergence. Thus when we talk about finding the critical value, we mean the critical value defined by the numerical method used, not the {\it exact} value from the PDEs. This will be explained in detail, using the Klein-Gordon equation in the paper. 

Second, although the singular potential term is used to mimic the effect of the local defect, such a model is not perfect. Associated parameters, such as the strength and the location of the singular potential terms, have a great degree of uncertainty. Thus, determining the critical value needs to be understood within the context of the uncertainty analysis. In this paper, as the first step toward thoroughly analyzing these two issues, we propose to use the generalized polynomial chaos (gPC) method \cite{Ghanem,LeMaitre,XiuBook} for finding the critical value of both the Klein-Gordon and sine-Gordon equations and the associated statistical quantities such as mean solutions around the critical value. The gPC method has been  successfully implemented to analyze the uncertainty quantification for problems such as fluid flows with uncertain initial and  boundary conditions, sensitivity analysis, steady-state problems, etc. \cite{Xiu2002-1, Xiu2003, Xiu2004, Xiu2002-3}. 

In our proposed method, we consider the associated parameter that yields the critical phenomenon of the system as a random variable, and assume that we will try to find a critical parameter value with the given numerical setting. For the Klein-Gordon equation, the associated parameter is the strength or amplitude of the singular potential term; for the sine-Gordon equation, the associated parameter is the initial velocity of the soliton solution in the presence of the singular potential term that has a fixed amplitude. The original equations are then parameterized by the random variable and the solution of the parameterized equations is obtained by the gPC method. That is, the solution is expanded by the orthogonal polynomials in the random space. The associated orthogonal polynomials are determined by the probabilistic nature of the random variable, which is the critical parameter in our case. For example, if the random variable has the uniform probability density function (PDF), the associated polynomial is the Legendre polynomial; if the PDF is normal distribution, the associated polynomial is the Hermite polynomial \cite{LeMaitre,XiuBook}. 

We first consider the Klein-Gordon equation with the singular potential and find some conditions for the critical phenomenon. The Klein-Gordon equation with the singular source term has been considered in various studies. For example, in \cite{Kazama}, the singular source term is considered in adopting the model, the meson theory, that the particle mass varies with space and time. In \cite{Dunne}, two $\delta$-functions are used as the singular potentials for Schr\"{o}dinger equation and the corresponding eigenfunctions are derived using the Weyl-Titchmarsh-Kodaira (WTK) spectral theorem. In \cite{GH2004}, the Klein-Gordon equation with the singular defect was considered and the chaotic behavior of the solution is studied. However, not every singular potential term yields the critical phenomenon. We consider the self-interacting singular potential term and show the existence of the critical phenomenon with a certain type of boundary condition. That is, the energy norm increases, vanishes, or decreases if the parameter value exceeds, is equal to, or is less than the critical value. The critical value of the parameter is obtained analytically in this case, which makes it possible to do an error analysis. For the critical value, we assume the uniform PDF and use the Legendre polynomial for the expansion. The resulting deterministic equations are solved using the Chebyshev collocation method with the consistent formulation for the Dirac $\delta$-function \cite{Jung2009}. Once the expansion coefficients are obtained, the solution is reconstructed for various values of the parameter. Using the fact that the energy norm remains constant for the critical value when the time is large, the invariant point is found as the critical value. The mean solution is also found by the first gPC mode. We compare the gPC approach and the MC approach. 

Then we apply the same method for the sine-Gordon equation. The sine-Gordon equation with singular potential terms is also found in much of the existing literature. The singularly perturbed solution of the sine-Gordon equation shows many interesting mathematical structures, such as the fractality, self-similarity, etc. \cite{GH2004,GH2007}. Moreover, such an equation yields the critical phenomenon associated with the critical parameter, as explained above briefly. In our work, we use the initial velocity of the soliton solution as the random variable, and try to find the critical initial velocity around which the particle-pass, particle-capture, and particle-reflection occur. For the random variable, we consider both the Hermite and Legendre polynomials--that is, the normal and the uniform distributions as the PDF. Unlike the Klein-Gordon equation, the reconstruction of the solution using the gPC method is not obtained directly. This is because the receiving soliton signal has two extreme cases, i.e. \textit{transmitted} or \textit{not-transmitted}, and because of the nonlinearity. The extreme solution behavior renders the PDF as the two $\delta$-function-like PDFs. For the particle-capture, the amplitude of the tail at the domain exit is $2\pi$; it vanishes eventually for the particle-pass. Such strong discontinuity in the solution behavior causes the Gibbs phenomenon in the reconstruction. These two extreme cases are not separable in the reconstruction. In the paper, we show such non-separable reconstruction for both the Hermite and Legendre cases. Thus instead of using the reconstruction by the gPC method, we use the gPC mean for the determination of the critical value. In this case, the Legendre gPC method yields much faster convergence than the Hermite gPC method. Numerical results show that the Legendre gPC method determines the critical value very efficiently and accurately. The determined critical value agrees well with the value found by the many direct MC simulations. 

To the best of our knowledge, gPC analysis for the singularly perturbed PDEs by the Dirac $\delta$-functions has not been thoroughly studied in existing literature. In particular, for the nonlinear optics equations with a singular defect in the disordered media, the uncertainty quantification is necessary. This is because the defect is defined in the highly localized regime and thus contains a high degree of intrinsic uncertainty. In our previous work \cite{JungSong2011}, some linear wave equations with the singular source term were considered with the Legendre PC method. Our current work takes the first step to investigate the uncertainty quantification of nonlinear optics equations, with the source term confined in the highly localized regime. 

The paper is composed of the following sections. In Section 2, the Klein-Gordon equation with the singular source term is considered. In this section, we show that there exists a critical phenomenon and find the exact value of the critical parameter value. In Section 3, we briefly explain the gPC method and derive the gPC expansion of the Klein-Gordon equation using the Legendre polynomials. In Section 4, numerical methods used for the numerical approximation of the deterministic equations derived in Section 3 are explained including the consistent method for the approximation of the Dirac $\delta$-function. In Section 5, numerical results of the gPC method for the Klein-Gordon equation are given. Comparison between the MC method and the gPC method is given. In Section 6, the critical phenomenon for the sine-Gordon equation is explained. The method of finding the critical initial velocity is explained using the PDF of the receiving signals and the gPC mean value. The Hermite gPC method is presented. In Section 7, the Legendre gPC is used for the sine-Gordon equation and the deterministic equations are derived.  Numerical results and remark on the error are presented. In Section 8, we briefly outlined the gPC method for the case that the randomess is in the amplitue of the singular potential function and explain that the same gPC method can be applied to this case. In Section 9, concluding remark and future works are presented. 

\section{Critical phenomenon in linear equation with a singular potential term} 
\subsection{Inhomogeneous Klein-Gordon equation}
We first consider the one-dimensional, inhomogeneous, Klein-Gordon type equation in the  general form given by 
\begin{eqnarray}
\label{inhom_KG}
u_{tt}-c^2u_{xx} + d^2u = p(x,t), \quad x \in \mathbb{R}, \quad t\in \mathbb{R}^{+}, 
\end{eqnarray}
where $u, p: \mathbb{R}\times \mathbb{R}^{+} \rightarrow \mathbb{R}$,  $c, d \in \mathbb{R} $. The initial conditions are  
\begin{eqnarray}
u(x,t) &=& u^0(x),  \quad t = 0, \nonumber \\
u_{t}(x,t) &=& 0, \quad t = 0, \nonumber   
\end{eqnarray}
where $u^0, u \in L^2(\mathbb{R})$ 
\begin{eqnarray}
u(x,t) \rightarrow 0,  \quad \quad |x| \rightarrow +\infty, \quad t\ge 0. \nonumber
\end{eqnarray}
The solution of Eq. (\ref{inhom_KG}) is given by
\begin{eqnarray}
u(x,t) =\frac{1}{2c} \int_{0}^{t}d\tau\int_{x-c(t-\tau)}^{x+c(t-\tau)}p(\xi,\tau)J_{0}\left[\frac{d}{c}\left[c^2(t-\tau)^2-(x-\xi)^2\right]\right] d\xi,
\end{eqnarray}
where $J_{0}$ is the Bessel function of the first kind.
The above solution can be derived using the Green's function method \cite{Evans}.

Now suppose that the source term involves the singular function such as $p(\xi,\tau) = \eta\delta(\xi)$ where $\eta \in \mathbb{R}$ is a constant and 
$\delta(\xi)$ is the Dirac $\delta$-function.   Then the solution becomes
\begin{eqnarray}
u(x,t) =\frac{\eta}{2c} \left[ \int_{0}^{t}J_{0}\left[\frac{d}{c}\left[c^2(t-\tau)^2-x^2\right]\right] d\tau\right]H(c(t-\tau)-x) , 
\end{eqnarray}
where $H$ is the Heaviside function and $c \ne 0$. 
As shown in the solution form, the parameter $\eta$ only scales the solution but does not yield the critical behavior of the solution. 
In order to obtain the critical solution, we consider the following several different cases. 

\subsection{Steady-State Solution: $c = 1, d = 1, p= \eta\delta(x)$}
We first consider the case of a possible steady-state solution with a singular source term.
Consider the case that $c = 1, d = 1, p = \eta \delta(x)$ for the bounded domain $x \in  \Omega = [-1,1]$ with the boundary 
conditions
$$
    u(-1,t) = 0  = u_x(1,t). 
$$
Then there exists a steady-state solution that satisfies 
\begin{eqnarray}
\label{steady_state1}
-u_{xx} + u = \eta\delta(x), \quad u(-1,\cdot)=0, \quad u_{x}(1,\cdot)=0. \nonumber 
\end{eqnarray}
By using these boundary conditions, the general solution can be easily derived as
$$ u(x) = \eta \left[e\left(\frac{e^2+1}{e^4+1} \right)\sinh (x+1) - H(x)\sinh(x) \right]. $$

If different boundary conditions are chosen such as 
$$
u(-1,\cdot)=1, \quad u_{x,\cdot}(1) = 0,
$$ 
then the solution becomes 
$$
u(x) = \left[\frac{\eta}{2}\frac{e^2(e^2+1)}{e^4+1} + \frac{e}{e^4+1} \right]e^{x} + \frac{1}{e^2}\left[ \frac{e^5}{e^4+1}-\frac{\eta}{2}\frac{e^2(e^2+1)}{e^4+1}\right]e^{-x} - \eta\sinh(x)H(x).  
$$
%

These examples show that there exists the steady-state solution for any value of $\eta$ but we do not 
observe the critical phenomenon in the solution induced by the singular potential term $p(x) = \eta \delta(x)$. 

\subsection{Self interacting potential: $c = 1, d = 0, p(x,t) = V(x)u(x,t)$ }
For another possible steady-state solution, we consider the self-interacting potential term given by 
\begin{equation} 
   u_{tt}-u_{xx}=Vu, \quad x \in \Omega, \quad t>0, 
   \label{self_potential}
\end{equation} 
 where $V : \mathbb{R} \rightarrow \mathbb{R}$ is the potential function. We show that some $V$ yield a critical phenomenon.

Similar equation has been used in quantum physics. When $V$ is given by the singular source term, the above equation represents some physical model such as one-dimensional diatomic molecule model with the singular potential well in the presence of static electric field \cite{Dunne}.

The boundary conditions we choose are 
\begin{eqnarray}
\left\{ \begin{array}{cl} u(x,t) = 0,  &  x=-1,\quad t>0, \\
u_{t}+u_{x} = 0,  &   x=1, \quad t>0. \end{array} \right.    
\label{self_potential_BC}
\end{eqnarray}
with some non-vanishing initial condition 
$$
  u(x,0) = u^0(x) \ne 0. 
$$
\noindent
The motivation of the given boundary condition of $(\partial_{t} + \partial_{x})u = 0$ at $x = 1$ is that if the potential term $V$ is localized inside $\Omega $ and vanishes at the boundary, then the given boundary condition yields the outflow boundary condition. 
The potential term $V$ will 
contain the singular term such that the potential effect vanishes at the domain boundary $x \in  \partial \Omega =\left\{ \pm 1 \right\}$.  

To show the criticality of the solution, we define the energy $E(t)$ of the system given by 
%
 $$
E(t) = \frac{1}{2}\left(\parallel u_{t} \parallel^{2}_{L^{2}} + \parallel u_{x} \parallel^{2}_{L^{2}}\right).
 $$ 
The rate of change of energy $E$ with respect to time $t$, $\dot{E}(t)$ is then simply computed by 
\begin{eqnarray} 
\dot{E}(t) = \int_{-1}^{1}u_{t}u_{tt}dx + [u_{t}u_{x}]_{-1}^{1}-\int_{-1}^{1}u_{t}u_{xx}dx.      \nonumber 
\end{eqnarray}
Using the boundary conditions we have 
\begin{eqnarray}
 \dot{E}(t) = \int_{-1}^{1} uu_{t}V dx -u_{t}^{2}(1,t). \nonumber  
 \end{eqnarray}
Suppose that the potential function $V$ is given by 
\begin{equation}
   V  = \eta \delta(x-\alpha) u(x,t), \quad \eta > 0, \quad \alpha \in (-1,\;1).
   \label{delta_potential} 
\end{equation}
Then 
\begin{equation}
    \dot{E}(t,\eta) = \eta \int_{-1}^1 uu_t \delta(x-\alpha) dx - u_t^2(1,t) = \eta u(\alpha,t)u_t(\alpha,t) - u_t^2(1,t). 
    \label{energy_change}
\end{equation}
Without loss of generality, we choose $\alpha = 0$. Using the continuity of $u$ at $x = \alpha = 0$, the jump condition is given by 
\begin{eqnarray}
   \lim_{\epsilon \rightarrow 0^+} 
   \left[ u_x(-\epsilon, t) - u_x(+\epsilon, t) \right]= \eta u(0,t), \quad \forall t \in \mathbb{R}^+. 
   \label{jump}
\end{eqnarray}

\begin{lemma}
The steady-state solution is determined by a unique $\eta^c$ and $\eta^{c} = \frac{1}{1+\alpha}$, $\alpha \in (-1,1)$.
\end{lemma}
\noindent
{\it Proof}: $\;$ Using the boundary conditions, we easily show that the steady-state solution $u^s(x)$ is given by 
\begin{eqnarray}
\label{ss-solution}
  u^s(x) = \left\{ \begin{array}{cc} C(x + 1) &  x \in [-1,\alpha] \\
                                                      C(1 + \alpha) & x \in [\alpha, 1] \end{array} \right. , 
\end{eqnarray}
where the constant $C$ depends on the initial condition. Notice that the parameter $\eta$ does not appear in the steady-state solution.  Using the jump condition, we obtain
$$
           \lim_{\epsilon \rightarrow 0^+} \left[ u^s_x(-\epsilon) - u^s_x(\epsilon) \right]= \eta u^s(\alpha). 
$$
Using Eq. (\ref{ss-solution}), the critical value of $\eta$, $\eta_c$ yielding the steady-state solution is obtained uniquely
$$
   \eta^c =  \lim_{\epsilon \rightarrow 0^+} \left[ u^s_x(-\epsilon) - u^s_x(\epsilon) \right]/ u^{s}(\alpha)= {1\over {1 + \alpha}}.\quad \Box
$$

\begin{theorem}
There exists a critical behavior around $\eta^c$ obtained in Lemma 1, that is, 
\begin{eqnarray}
\dot{E}(t,\eta) > 0 \quad \mathrm{if} \; \eta > \eta^{c}, \nonumber \\
\dot{E}(t,\eta) = 0 \quad \mathrm{if} \; \eta = \eta^{c}, \nonumber \\
\dot{E}(t,\eta) < 0 \quad \mathrm{if} \; \eta < \eta^{c}, \nonumber 
\end{eqnarray}
for  $t \rightarrow \infty$ where $\eta \in N_{\epsilon}(\eta^{c})$ for an arbitrary small value of $\epsilon > 0$.\\
\end{theorem}

\noindent
{\it Proof}: 
Without loss of generality we consider $\alpha = 0$, i.e.,  the position of $\delta$-function is at the origin and suppose that we choose 
$$
         \eta = \eta^c + \Delta \eta = 1 + \Delta \eta, 
$$
where $\eta^{c} = 1$ (from Lemma 1) and $0< |\Delta \eta| \ll 1$. Based on these assumptions, we look for the solution in terms of the power series in $\Delta \eta,$ that is, 
$$
     u = u^s + \Delta\eta u_1 + (\Delta\eta)^2 u_2 + (\Delta\eta)^3 u_3 + \cdots, 
$$
where $u^{s}$ is the steady-state solution. Plug $u$ into the PDE given by Eq. \ref{self_potential} with \ref{delta_potential} and then  each of $(\Delta \eta)^{n}$ terms yields
\begin{eqnarray}
 (\Delta\eta)^1: \quad     (u_1)_{tt} - (u_1)_{xx} &=& \delta(x)(u^s + u_1), \nonumber \\
    (\Delta\eta)^2: \quad        (u_2)_{tt} - (u_2)_{xx} &=& \delta(x)(u_1 + u_2), \nonumber \\
     (\Delta\eta)^3: \quad            (u_3)_{tt} - (u_3)_{xx} &=& \delta(x)(u_3 + u_2), \nonumber \\
             \cdots   \nonumber 
\end{eqnarray}
We choose the initial condition such that the steady-state solution becomes 
\begin{eqnarray}
  u^s(x) = \left\{ \begin{array}{cc} (x + 1) ,&  \mbox{ for } x \in [-1,0] \\
                                                      1,  & \mbox{ for } x \in [0, 1] \end{array} \right. . 
                                                      \nonumber
\end{eqnarray}
\noindent
For the PDE of the first order term above, the boundary conditions and the continuity condition yields 
\begin{eqnarray}
u_{1}(x,t) = \left\{ \begin{array}{cc} a(1 + x) + ct(1 + x) , & \mathrm{for} \; x \le 0 \\
               a  -cx + ct , &  \mathrm{for} \; x > 0  \end{array} \right. . \nonumber
\end{eqnarray}
The jump condition, $-u_{x}(\epsilon^{+}) + u_{x}(\epsilon -) = u^{s}(0) + u_{1}(0)$ yields $c=1$. We choose the initial condition such that $a = 1$. Then by ignoring the higher-order terms $(\Delta \eta)^{(n)}$, $n = 2, 3, 4, \cdots$, we obtain
%
\begin{eqnarray}
u(x,t) = \left\{ \begin{array}{cc} 1+ x + \Delta\eta \left( 1 + x + t(1 + x)\right) , & \mathrm{for} \; x \le 0, \\
             1+ \Delta \eta ( 1  -x + t) , &  \mathrm{for} \; x > 0.  \end{array} \right.  \nonumber
\end{eqnarray}
\noindent
If $\Delta \eta > 0$ the solution increases unboundedly and the energy increases with time. If $\Delta \eta = 0$, the solution becomes the steady-state solution and the energy is conserved. If $\Delta \eta <0$, the solution decays until $t = - \left(\frac{1}{\Delta\eta} + 1 \right) $. For $t \geq - \left(\frac{1}{\Delta\eta} + 1 \right)$, the higher-order term becomes significant and they should be considered so that the energy further decays. $\Box$

\begin{remark}
To get the complete solution, one needs to solve the higher-order terms. 
It may become difficult to solve those PDEs exactly but the proof is enough to show that there exists a critical behavior of the solution around $\eta^c$  when $|\Delta \eta | \ll 1$. 
For example, if $\Delta \eta \approx 5\times 10^{-5}$, the first-order term is dominant up to $t \approx 0.2 \times 10^5$, which is enough for the long time behavior of the solution for our purpose. 
\end{remark}

\begin{figure}
	\centering
		\includegraphics[width=0.45\textwidth]{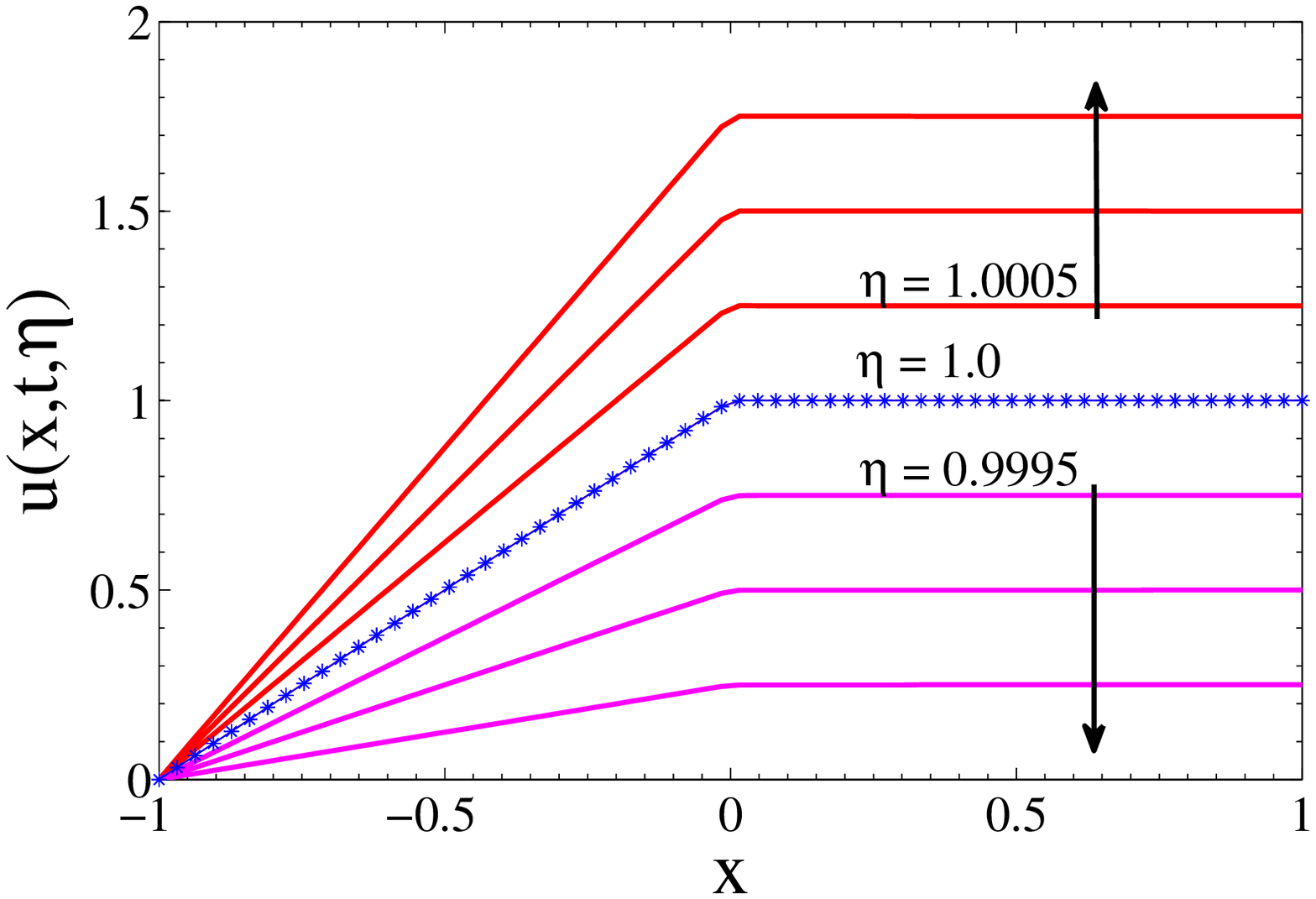}
		\includegraphics[width=0.43\textwidth]{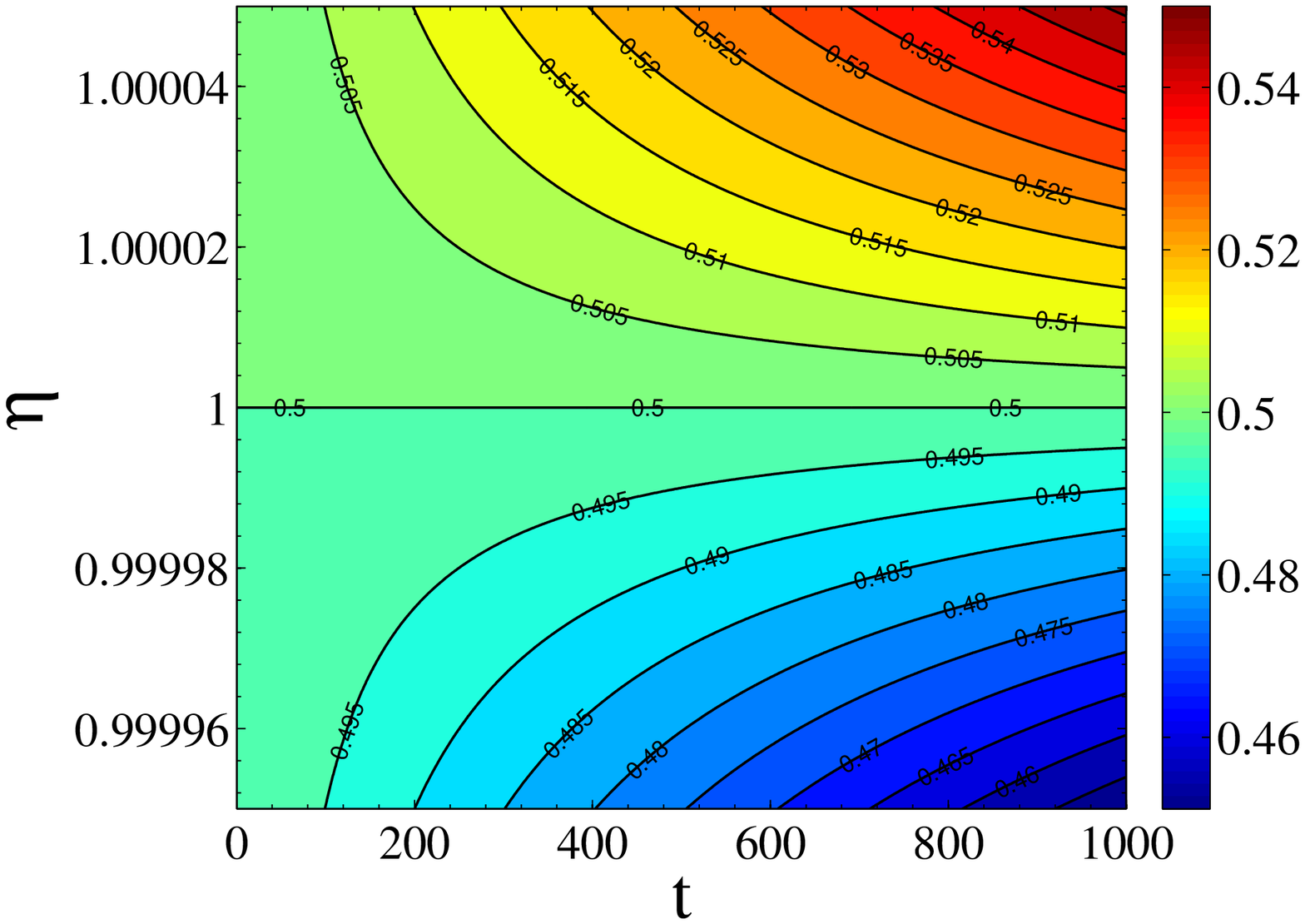}
		\caption{Left: Analytical solutions of $u(x,t)$ at different final times $t_{f} = 0,\; 500,\; 1000,\; 1500$.  Right: Contour plot of $E$ for $\eta \in [0.9995, \; 1.0005]$}
	 \label{exact_solution}
\end{figure}

The energy of the system up to the first order term is given by 
\begin{eqnarray}
E(t,\eta) &=& \frac{1}{2}\left[\eta + t(\eta-1)\right]^2 + \frac{7}{6}\left(\eta -1\right)^{2} + (1+t)(1-\eta)^{2}\left[\eta + t(\eta -1)\right]. \nonumber 
\end{eqnarray}
And the rate of change of $E$ with time is then given by 
$$
{\dot{E}} = \eta\left(\eta -1\right)\left[ \eta + t \left( \eta -1\right) \right]+ \left(\eta-1\right)^{3}\left(1+t\right) .
$$
Figure \ref{exact_solution} shows the solutions at different times $t_{f} = 0,\; 500,\; 1000,\; 1500$ for different values of $\eta$,  $\;\eta = 0.9995,\; 1.0005$ around $\eta^c = 1$ (left) and the contour of $ E$ in $\eta$ and $t$ plane (right). 

\subsubsection{Linearized sine-Gordon equation}
The same phenomenon can be obtained for the linearized sine-Gordon equation with the point-like source term where $\sin(u) \approx u$ such that 
\begin{eqnarray}
\label{linear_SG}
u_{tt}-u_{xx}=(\eta\delta(x)-1)u,  \quad x \in (-1,1), \quad t > 0,   
\end{eqnarray}
with the same boundary conditions 
\begin{eqnarray}
u(-1,t) = 0, \quad \mbox{ and } \quad u_{t}+u_{x}=0 \quad \mbox{ at }  x = 1.\nonumber
\end{eqnarray}
%

As in above, the energy, 
$
E(t) =   \frac{1}{2}\int_{-1}^{1}\left(u_{t}^2+u_{x}^2 \right)dx, 
$
has the rate of change below 
\begin{eqnarray}
\label{sg_energychange}
\dot{E}(t) 
&=&\int_{-1}^{1}-uu_{t}(x,t) dx + \eta u_t(0,t)u(0,t)-u_{t}^2(1,t).  \end{eqnarray}
We have a similar phenomenon 
$
\dot{E}(t)\lesseqqgtr 0 \Rightarrow \int_{-1}^{1}-uu_{t}(x,t)  dx +\eta  u(0,t) u_{t}(0,t)-u_{t}^{2}(1,t)\lesseqqgtr 0, 
$
that occurs around the following steady-state solution 
\begin{eqnarray}
u^{s}(x) = \left\{\begin{array}{cc} C(1+e^{2})(e^{-x} - e^{x+2}),&\mathrm{for} \quad x<0, \\
   C(1-e^{2})(e^{x} + e^{2-x}), & \mathrm{for} \quad x>0, \end{array} \right. 
\end{eqnarray}
where $C$ is a constant that depends on the initial condition. The critical value of $\eta$ 
is $\eta^c = 2\coth(2)$.

\section{Expansion using the gPC method}
\subsection{gPC method and Askey scheme.}


 We use the gPC method for the solution of the Klein-Gordon and sine-Gordon equations using the Wiener-Askey scheme \cite{Xiu2002-2}, in which Hermite, Legendre, Laguerre, Jacobi, generalized Laguerre orthogonal polynomials are used for modeling the effect of continuous random variables described by normal, uniform, exponential, beta and gamma probability distributions,  respectively \cite{XiuBook,Xiu}.  These orthogonal polynomials are optimal for those PDFs since the weight function in the inner product and its support range correspond to the probability density functions for those continuous distributions.

Table (1) shows the set of polynomials which provide an optimal basis for different continuous PDFs \cite{LeMaitre,XiuBook}. It is derived from the family of hypergeometric orthogonal polynomials known as the Askey scheme, for which the Hermite polynomials originally employed by Wiener \cite{Wiener} are a subset. 

{{
\begin{table}
\label{table1}
\caption{Correspondence between standard forms of continuous probability distributions and Askey scheme of continuous hypergeometric polynomials \cite{LeMaitre,XiuBook}.}
\begin{tabular}
{ccccc}\hline\em Distribution &\em  Density Function  &\em  Polynomial &\em Weight function &\em Support range \\\hline\hline Normal & $\frac{1}{\sqrt{2\pi}}e^{-\frac{x^2}{2}}$  & Hermite $H_{n}(x)$ & $e^{-\frac{x^2}{2}}$ & $(-\infty \; \infty)$  \\\hline Uniform & $\frac{1}{2}$ & Legendre $P_{n}(x)$ & $1 $ & $[-1 \; 1] $  \\\hline Beta & $\frac{(1+x)^{\alpha}(1+x)^{\beta}}{2^{1+\alpha +\beta}B(\alpha+1,\beta+1)}$ & Jacobi $P_{n}(\alpha, \beta)(x)$ & $(1+x)^{\alpha}(1+x)^{\beta} $ & $[-1 \; 1] $  \\\hline Exponential & $e^{-x}$ & Laguerre $L_{n}(x)$ & $e^{-x}$ & $[0 \; \infty)$ \\\hline Gamma & $\frac{x^{\alpha}e^{-x}}{\Gamma(\alpha + 1)}$ & Generalized Laguerre $L_{n}^{(\alpha)}(x)$ & $x^{\alpha}e^{-x}$ & $[0 \; \infty)$  \\\hline  
\end{tabular}
\end{table}
}}

Following the standard GPC expansion, we assume that $u(x,t, \xi)$ is sufficiently smooth  in $\xi$ and has a converging expansion of the form 
$$(u,x,\xi) = \sum_{k=0}^{\infty}\hat{u}_{k}(x,t)P_{k}(\xi),$$
where the orthonormal polynomials $P_{k}(\xi)$ correspond to the PDF of the random variable $\xi$ and satisfy the following orthogonality relation:
$$
\mathbf{E}[P_{k}P_{l}] := \int P_{k}(\xi)P_{l}(\xi)\rho(\xi)d\xi = \delta_{kl}, 
$$
where $\delta_{kl}$ is the Kronecker delta and $\rho(\xi)$ is the weight function. Note that the polynomials are normalized. 



\subsection{Expansion using the Legendre gPC for the Klein-Gordon equation}

We use the gPC method to find the critical value of $\eta$ numerically. 
Assume that $\eta > 0$ is a random variable and $\eta \in \left[ a ,\; b \right]$ has the uniform distribution.
Let $\xi$ be a random variable with the uniform PDF over $\Omega = [-1,1]$. Let $\eta$ be a linear map, $ \eta: \Omega\rightarrow [a,\;b] $, $\eta =  \frac{b-a}{2}\xi + \frac{a+b}{2}$.

The solution $u(x,t,\xi)$ will be expanded using the Legendre polynomial in $\xi$ with time $t$ and the value of $\xi$ will be sought that yields the invariant solution with time when $t\rightarrow \infty$. Such $\xi$ or $\eta$ is then determined as the critical value for the Klein-Gordon equation. 

For this, we first expand the solution $u$ in the Legendre polynomial in the Galerkin sense such as in \cite{GottliebXiu}:
%
\begin{equation} 
u(x,t,\xi)= \sum_{l=0}^{\infty}\hat{u}_{l}(x,t)L_{l}(\xi), 
\label{eq33}
\end{equation}
where $\hat{u}_{l}(x,t)$ are the expansion coefficients and $L_l(\xi)$ the Legendre polynomials \cite{Bateman}. 
We consider the truncated version of Eq. (\ref{eq33}) by considering the first $N+1$ terms. 
Plugging Eq. (\ref{eq33}) into Eq. (\ref{self_potential}) and using $V$ from Eq. (\ref{delta_potential}) with $\alpha = 0$, we get 
$$\sum_{l=0}^{N}\left[\left\lbrace\hat{u}_{l}(x,t) \right\rbrace_{tt} - \left\lbrace \hat{u}_{l}(x,t) \right\rbrace_{xx}  \right ]L_{l}(\xi) = \left(\frac{b-a}{2}\xi + \frac{a+b}{2}\right)\delta(x)\sum_{l=0}^{N}\hat{u}_{l}(x,t)L_{l}(\xi). $$
Multiplying both sides by $L_{l^{'}}(\xi)$ yields

$$
\sum_{l=0}^{N}\left[\left\lbrace\hat{u}_{l}(x,t) \right\rbrace_{tt} - \left\lbrace \hat{u}_{l}(x,t) \right\rbrace_{xx}  \right ]L_{l}(\xi)L_{l^{'}}(\xi) = \left(\frac{b-a}{2}\xi + \frac{a+b}{2}\right)\delta(x)\sum_{l=0}^{N}\hat{u}_{l}(x,t)L_{l}(\xi)L_{l^{'}}(\xi).
$$
%
Integrating both sides over $[-1,1]$ and using the orthogonal property of the Legendre polynomials, we get,
\begin{eqnarray}
\label{gPC_legendre1}
&& \left[\left\lbrace\hat{u}_{l}(x,t) \right\rbrace_{tt} - \left\lbrace \hat{u}_{l}(x,t) \right\rbrace_{xx}  \right ]\frac{2}{2l+1} 
 = \frac{2}{2l+1}\left(\frac{a+b}{2}\right)\delta(x)\hat{u}_{l}(x,t) \nonumber \\
 &+& \left(\frac{b-a}{2}\right)\delta(x)\sum_{l=0}^{N}\hat{u_{l}}(x,t)\int_{-1}^{1}\xi L_{l}(\xi)L_{l^{'}}(\xi)d\xi.
\end{eqnarray}
For the second term in the right hand side of the above equation, one can show the following relation
$$
\int_{-1}^{1}\xi L_{l}(\xi)L_{l^{'}}(\xi) d\xi= \left\lbrace \begin{array}{ccc}\frac{2(l+1)}{(2l+1)(2l+3)} & \quad \mathrm{for} & l = l^{'}+1 \\ \frac{2l}{(2l+1)(2l-1)} & \quad \mathrm{for} & l = l^{'}-1 \\ 0 & \mathrm{otherwise}\end{array} \right. 
\label{legendre_integration}
$$
(for the derivation of the relation, see Appendix 2). 

Then, Eq. (\ref{gPC_legendre1}) can be rewritten as 
\begin{eqnarray}
\left\lbrace \hat{u}_{l}(x,t)\right\rbrace _{tt} &-& \left\lbrace \hat{u}_{l}(x,t)\right\rbrace _{xx} = \delta(x)\left[\left(\frac{b-a}{2}\right) \frac{l}{2l-1}\hat{u}_{l-1}(x,t)  \right.
\nonumber \\
&+& \left. \left(\frac{a+b}{2}\right)\hat{u}_{l}(x,t) + \left(\frac{b-a}{2}\right) \left(\frac{l+1}{2l+3}\right)\hat{u}_{l+1}(x,t)\right]. 
\label{gPC_legendre2}
\end{eqnarray}
Substitute Eq. (\ref{eq33}) into the left boundary condition of Eq. (\ref{self_potential_BC})  to  get, 
$$
\sum_{l=0}^{N}\hat{u}_{l}(-1,t)L_{l}(\xi)=0. 
$$
By using the  orthogonal property of the Legendre polynomials, we have 
\begin{equation}
\hat{u}_{l}(-1,t)=0, \quad \forall l = 0, \cdots, N. 
\label{gpc_legendre_bcl}
\end{equation}
Similarly the right boundary condition of Eq. (\ref{self_potential_BC}) yields
$$
\sum_{l=0}^{N}\left[\left\lbrace \hat{u}_{l}(1,t)\right\rbrace _{t} + \left\lbrace \hat{u}_{l}(1,t)\right\rbrace _{x}\right]L_{l}(\xi)=0.
$$
By using orthogonality we have
\begin{equation}
\left\lbrace \hat{u}_{l}(1,t)\right\rbrace _{t} + \left\lbrace \hat{u}_{l}(1,t)\right\rbrace _{x} = 0, \quad \forall l = 0, \cdots, N. 
\label{gpc_legendre_bc2}
\end{equation}
We choose the initial condition as 
$$(1+x) = \sum_{l=0}^{N}\hat{u}_{l}(x,0)L_{l}(\xi).$$ 
Again using the orthogonal property, we get 
$$\left(1+x \right)\int_{-1}^{1}L_{l}(\xi)d\xi = \frac{2}{2l+1}\hat{u}_{l}(x,0). $$
This can be written as 
\begin{eqnarray} \hat{u}_{0}(x,0)&=&(1+x), \nonumber \\
                \hat{u}_{l}(x,0) &=& 0, \quad \forall l = 0, \cdots, N. 
                \label{gpc_legendre_initial}\end{eqnarray}
So from Eqs. (\ref{gPC_legendre2}), (\ref{gpc_legendre_bcl}), (\ref{gpc_legendre_bc2}), and (\ref{gpc_legendre_initial}) we have  $(N+1)$ equations for the gPC solution of the Klein-Gordon equation. 
%
%
%

\section{Numerical methods for solving PDEs with singular source terms}
\subsection{Consistent Chebyshev method}
The derived gPC equations for ${\hat u}_l(x,t)$ in Section 3 will be solved numerically. For the space derivative, the Chebyshev spectral collocation method is used based on the Gauss-Lobatto quadrature points $x_j  = -\cos(j\pi/m), j = 0, 1, 2, \cdots, m$ where $(m+1)$ is the total number of the quadrature points for $x_j \in [-1,1]$. As shown in Section 3 the derived gPC equations for ${\hat u}_l(x,t)$ still contains the $\delta$-function. Solving the PDE containing the singular source term properly is the crucial part for obtaining the accurate and stable numerical solutions for ${\hat u}_l(x,t)$. To deal with this problem, the $\delta$-function is approximated by the consistent formulation proposed in \cite{Jung2009} where the consistent formulation was proposed for the one-dimensional problems and it was shown that the method yields highly accurate results at the collocation points. It is not straightforward to extend such a formulation to higher dimensional problem, but as shown in the following sections and also shown in \cite{JungDon2009,JungSong2011,JungKhannaNagle2009}, this formulation yields an accurate result for some PDEs and can be expanded to higher dimension by using the dimension-by-dimension approach. Moreover, the method is simple and easy to implement. The underline motivation is that the approximation of the $\delta$-function is sought such that the steady-state solution can be obtained. This method is also referred as the well-balanced equation in the literature. The steady-state solution of the equation, $u^s(x)$ in the previous section satisfies
\begin{eqnarray}
           - u^s_{xx} = \eta \delta u^s(x). 
\label{steady}
\end{eqnarray}

Let $\delta_m$ be the approximation of $\delta$ on the $(m+1)$ collocation points.  Then the consistent formulation approximates $\delta_m$ as the derivative of the Heaviside function $H(x)$,  using the same spectral derivative operator. If $H_m$ denotes the Heaviside function defined at the collocation points $x_j$, then $\delta_m$ is given by the following: 
\begin{eqnarray}
        \delta(x) = {{d}\over{dx}}H \longrightarrow \delta_m := {\bf D} H_m, 
\end{eqnarray}
where ${\bf D} \in \mathbb{R}^{(m+1)\times (m+1)}$ is the Chebyshev derivative matrix and $H_m \in \mathbb{R}^{(m+1)\times 1}$ given by $H_m = (0, 0, \cdots, 0, 1, 1, \cdots, 1)^T$. To avoid any ambiguity, we choose the odd value of $m$ such that for $k = (m-1)/2$, 
$$
      \xi_k < 0 < \xi_{k+1}, \; \mathrm{and} \; \left(H_{m}\right)_{i} = \left\{\begin{array}{cc} 0 ,&\mathrm{for} \quad i \le k \\
   1, & \mathrm{for} \quad i \ge k+1 \end{array} \right. .
$$
That is, the location of the $\delta$-function is off the collocation point. Then the consistent formulation of Eq. (\ref{steady}) becomes
\begin{eqnarray}
\label{linear}
    -{\bf D}^2 U = \eta \Delta U, 
\end{eqnarray}
where $\Delta$ is the $(m+1)\times (m+1)$ matrix whose diagonal elements are $\delta_m$ and $U = (U(x_0), \cdots, U(x_m))^T$. 

\subsection{Numerical errors of numerically obtained critical values of $\eta^c$} 
\begin{figure}
	\centering
		\includegraphics[width=0.47\textwidth]{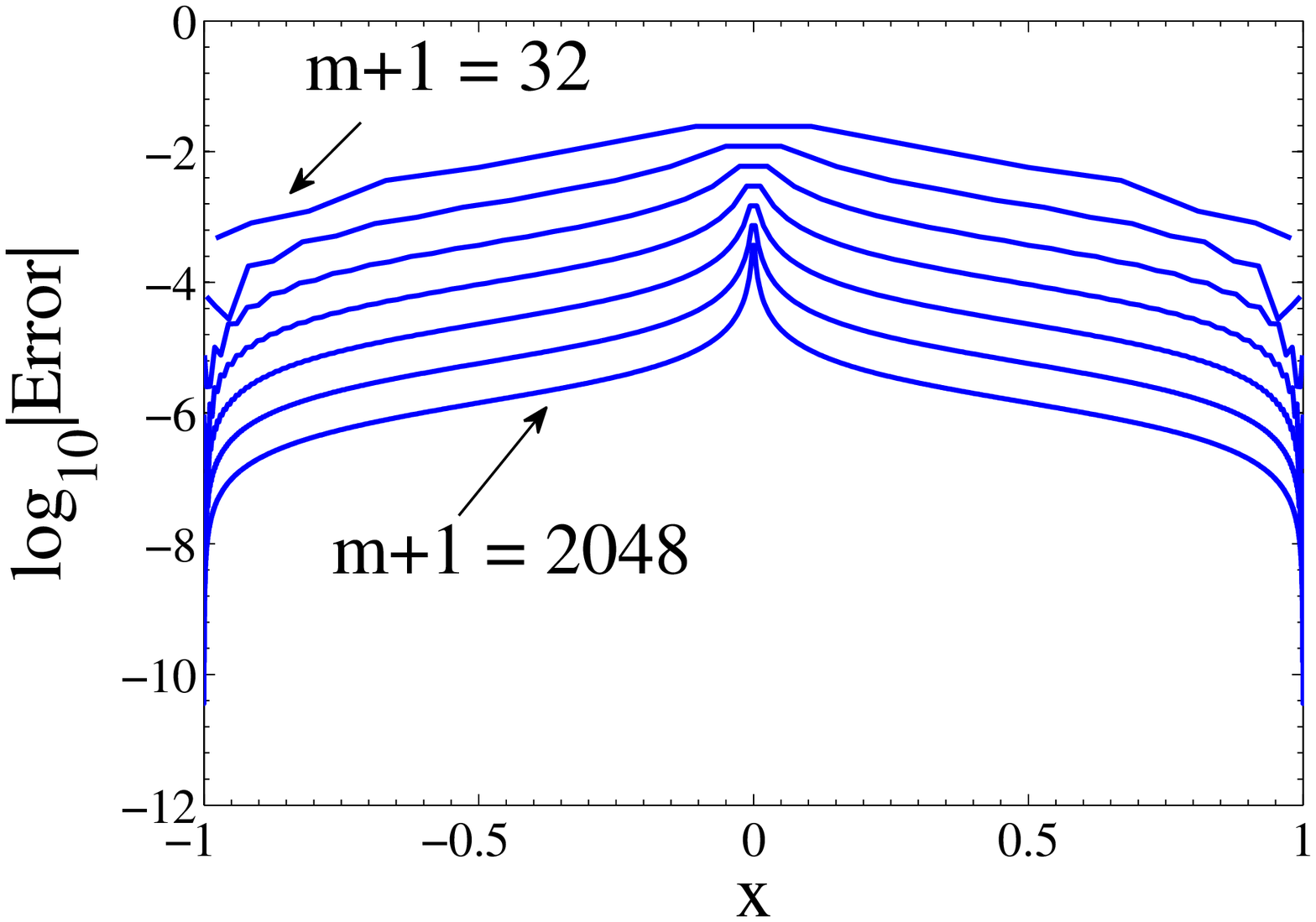}
		\includegraphics[width=0.47\textwidth]{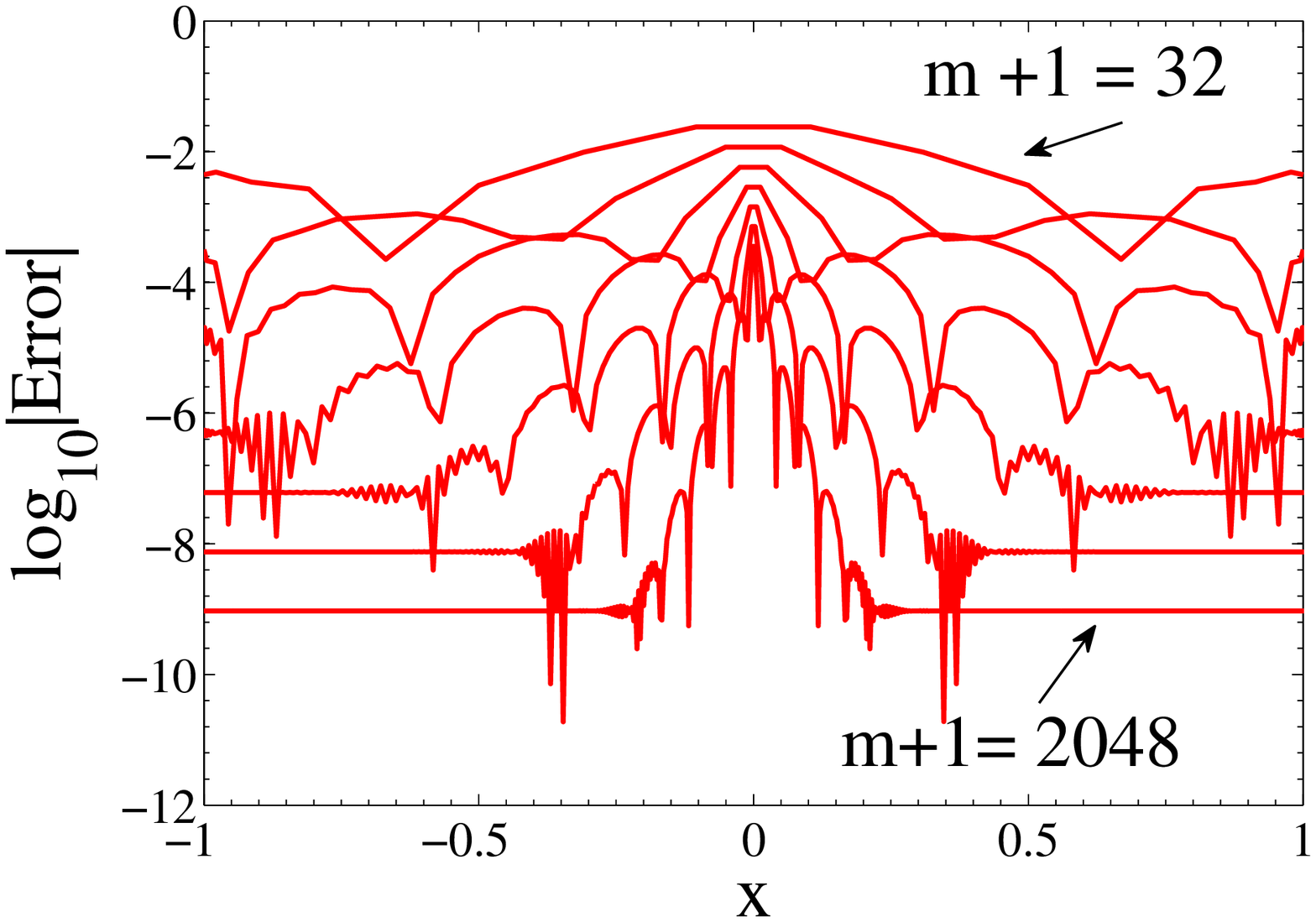}
		\caption{Pointwise errors of $|u-u^{e}|$. Left: Without the spectral filter. Right: With the spectral filter. $m = 32,\; 64\; \cdots, 2048$ from top to bottom.}
	 \label{fig:steady_state}
\end{figure}

As shown in Section 2, the critical value of $\eta^c$ is determined by the following 
$$
\eta^c =  \lim_{\epsilon \rightarrow 0^+} \left[ u^s_x(-\epsilon) - u^s_x(\epsilon) \right]/ u^{s}(\alpha). 
$$
That is, the critical value is determined by the local values of $u$ and $u_x$ near $x = 0$. Thus the convergence of the numerically computed value of $\eta^c$ depends on the local convergence of $u$ and $u_{x}$ near $x = 0$. Figure \ref{fig:steady_state} shows how the local convergence is obtained by our numerical approximation using the consistent method. In the figure, the convergence of the solution   is obtained by solving Eq. (\ref{linear}) iteratively. The left figure shows the pointwise errors with different values of $m$ and the right figure  the pointwise errors after the spectral filter is applied. As shown in these figures, the convergence can be enhanced away from the discontinuity when the filter operation is applied. But the overall convergence is slow due to the singularity at $x = 0$, which is the Gibbs phenomenon for the spectral method \cite{Gottlieb,Hesthaven}. This implies the slow convergence of the numerically obtained $\eta^c$.    

%
%
%

\begin{figure}
	\centering
		\includegraphics[width=0.43\textwidth]{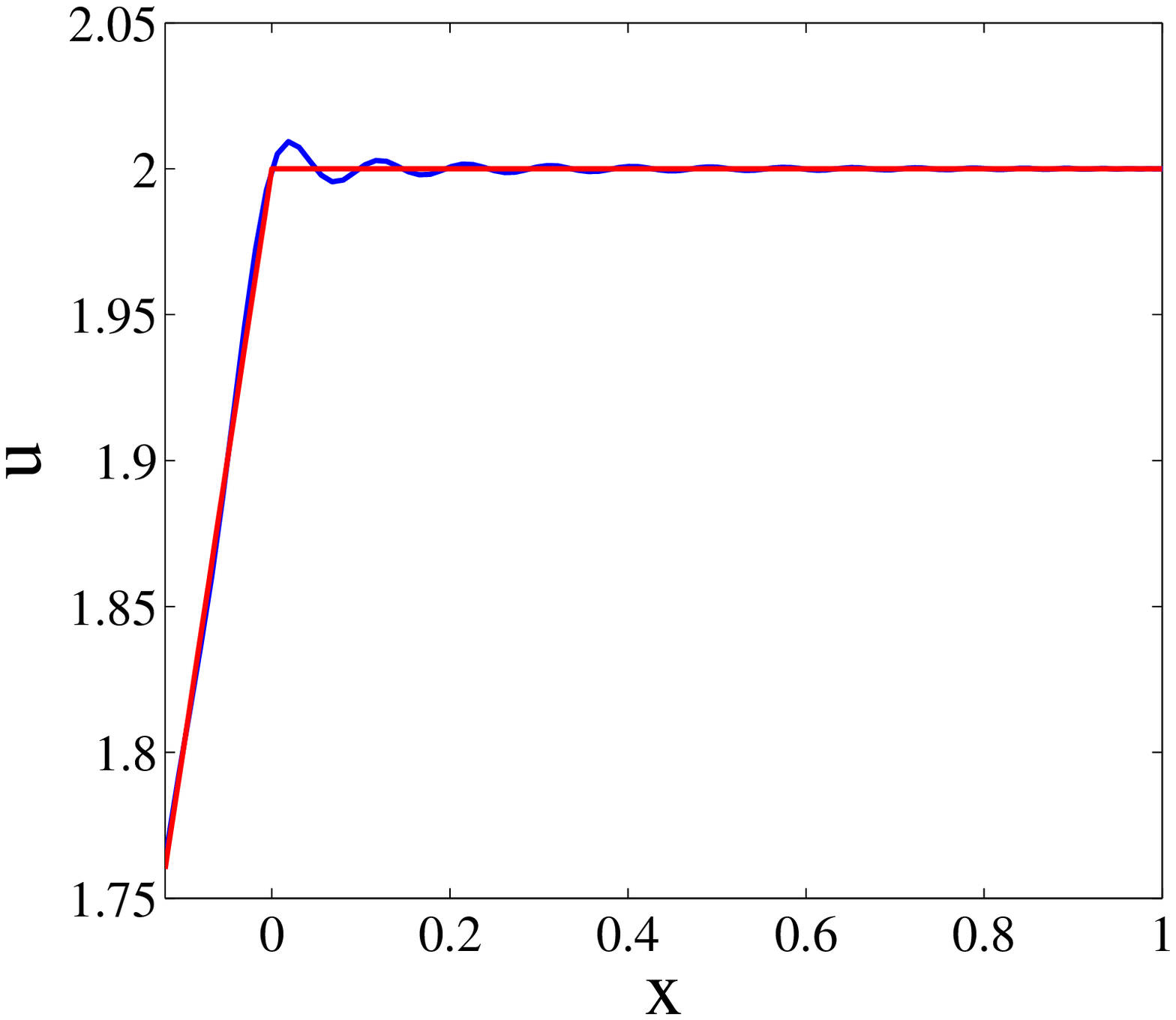}
		\includegraphics[width=0.43\textwidth]{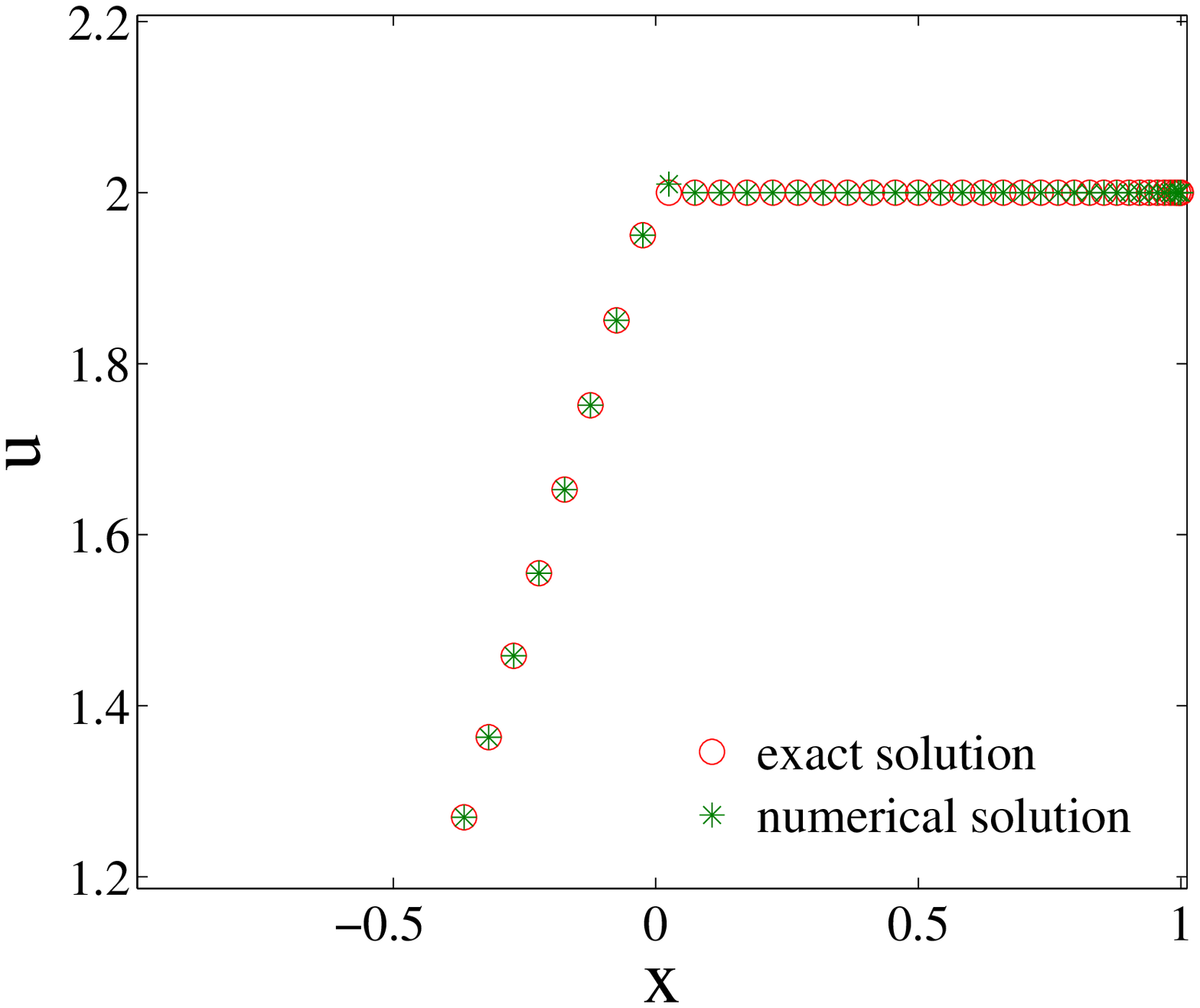}
		\caption{Left: The Chebyshev Galerkin solution.  Right: The Chebyshev collocation solution.}
	 \label{fig:cheb_leg}
\end{figure}

To show such Gibbs phenomenon we consider the exact solution  of Eq. (\ref{steady})  given by 
\begin{equation}
u^e(x) = A\left((x+1) - \int_{-1}^{x}H(\xi)d\xi\right).
\label{steady_state_sol}
\end{equation}
The solution is continuous in $(-1,1)$ but at $x=0$ it is not differentiable, $u(x) \in \mathrm{C}^{0}$. 
%
%
Let the exact solution $u^e(x)$ be expanded as the Gegenbauer polynomial of degree $m$: 
$$u^{e}_{m} = \sum_{k = 0}^{m}\hat{u}^{\mu}_{k}C^{\mu}_{k}(x),$$
where $$\hat{u}^{\mu}_{k} = \frac{1}{h^{\mu}_{k}}\int_{-1}^{1}(1-x^2)^{\mu - \frac{1}{2}}C^{\mu}_{k}(x)u^e(x)dx, \qquad 0 \leq k \leq m,$$
We use the particular case of the Gegenbauer polynomials for $\mu = 0$, which 
corresponds to the Chebyshev polynomials.
$\;$ Figure \ref{fig:cheb_leg} shows the Chebyshev Galerkin solution (left) and the Chebyshev collocation solution (right) with $m+1 = 128$.  As shown in these figures, there exist the Gibbs oscillations near $x = 0$. The Gibbs oscillation is seen mainly at the collocation point $x = x_{(m+1)/2}$ for the collocation solution. These oscillations on the collocation points yield the error in determining the critical value $\eta^c$. 

Table \ref{table2}  shows the convergence of $\eta^n$ with $m$, where $\eta^n$ is the numerically determined critical value of $\eta^c$.  The critical value we will find from the following section using the given value of $m$ will be those values in the table. For example, if we choose $m = 63$, the numerically computed critical value, $\eta^n$ is $\eta^n = 1.004332$ and our objective is to find that value as our critical value. That is,  given a finite value of $m$, there exists the corresponding $\eta^{n}$ uniquely and when we  find $\eta^{c}$, we mean it by finding $\eta^{n}$ instead of $\eta^c$. 

\begin{table}
\caption{Convergence of the critical value of $\eta^n$ with $m$. }
\label{table2}
\begin{center}
\begin{tabular}
{cc}\hline\em $m+1$ &\em  $\eta^n$   \\\hline\hline $32$ & $1.009595$    \\\hline $64$ & $1.004332$   \\\hline $128$ & $1.002036$   \\\hline $256$ & $1.000979$   \\\hline  
\end{tabular}
\end{center}
\end{table}


\subsection{Time integration} 
We use the second order centered difference in time for the time integration. For example let $U^n \in \mathbb{R}^{(m+1)\times 1}$ be the solution vector defined at $t = n\Delta t$ with the time-step $\delta t$ such as $U^n = (U(x_0, n\delta t), U(x_1, n\delta t), \cdots, 
U(x_m,n\delta t))^T$ for $u_{tt}-u_{xx}=P$.  Then the scheme is given by 
$$
     U^{n+1} = 2U^n - U^{n-1} + \delta t^2 {\bf D}\cdot {\bf D} U^n + P^n, 
$$
where the element of the column vector $P^n\in \mathbb{R}^{(m+1)\times 1}$ is $(P^n)_j = \eta (\delta_m)_j U^n_j$ and we apply $\bf D$ twice to $U^n$ for its second derivative in $x$. This is a two-step method in time. $U^0$ is given by the initial condition and we use the first-order approximation for $U^1$.

%
%
%
%

Then the numerical solution of $\hat{u}_l \in \mathbb{R}^{(m+1)\times 1}$ for each $l = 0, \;1,\;\cdots,N$ is obtained, by 
\begin{eqnarray}\hat{u}_{l}^{n+1} = \left[\Delta(\delta t)^{2}A_{l}\hat{u}^n_{l-1} +\left\lbrace B\Delta(\delta t)^{2}+ D^2(\delta t)^2 + 2I\right\rbrace \hat{u}^{n}_{l} + \Delta(\delta t)^{2}C_{l}\hat{u}^n_{l+1} \right] - \hat{u}^{n-1}_{l}, \nonumber \\
   \end{eqnarray}
where $A_{l} =  \left(\frac{b-a}{2}\right)\frac{l}{2l-1}$, $B = \frac{a+b}{2}$, $C_{l} = \left(\frac{b-a}{2}\right)\frac{l+1}{2l+3}$, and $\Delta = {\bf D}\cdot H_m$. The boundary conditions are 
\begin{eqnarray} 
\hat{u}^{n+1}_{l}(-1,t) &=& 0, \nonumber \\ \hat{u}^{n+1}_{l}(1,t) &=& (1-\lambda)u^{n}_{l}(1,t) + \lambda \hat{u}^{n}_{l-1}(1,t), \nonumber 
\end{eqnarray}
where $\lambda = \frac{\delta t}{\delta x}$ and $\delta x = x_m - x_{m-1}$.  For the boundary condition at $x = 1$, we use the first order finite difference both in time and space. 

\section{Numerical results: determination of $\eta_{c}$ and the mean solution}
\subsection{Critical behavior}
For the numerical example, we consider $(m+1)= 64$ and the corresponding $\eta \approx 1.004332$ from Table 2.

\begin{figure}
	\centering
		\includegraphics[width=0.45\textwidth]{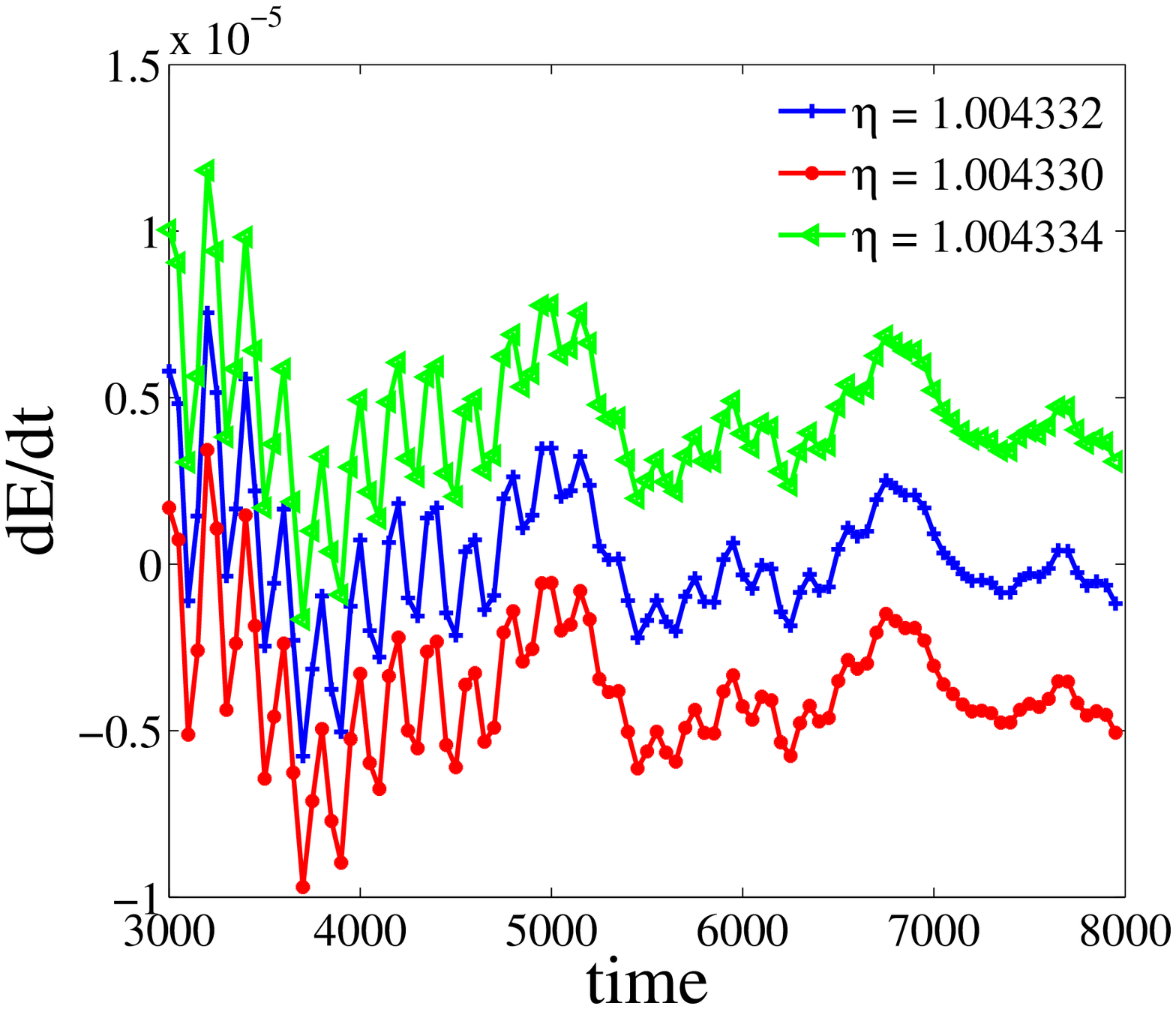}
		\includegraphics[width=0.50\textwidth]{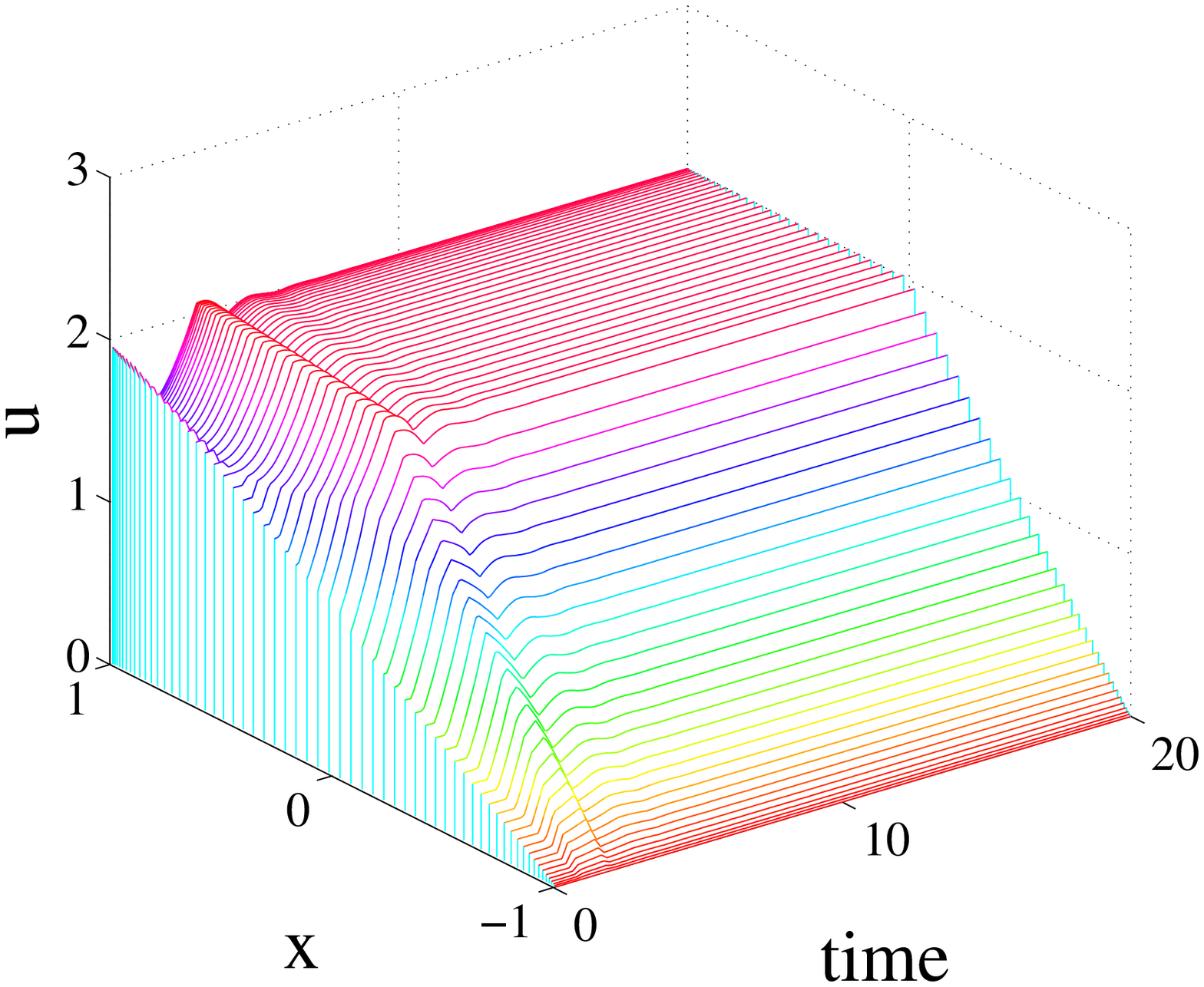}
		\caption{Left: $\dot{E}$ with different values of parameter $\eta$. Right: The steady-state solution for $\eta = 1.004332$
		}
	 \label{fig:Eta}
\end{figure}

\begin{figure}
	\centering
		\includegraphics[scale=0.2]{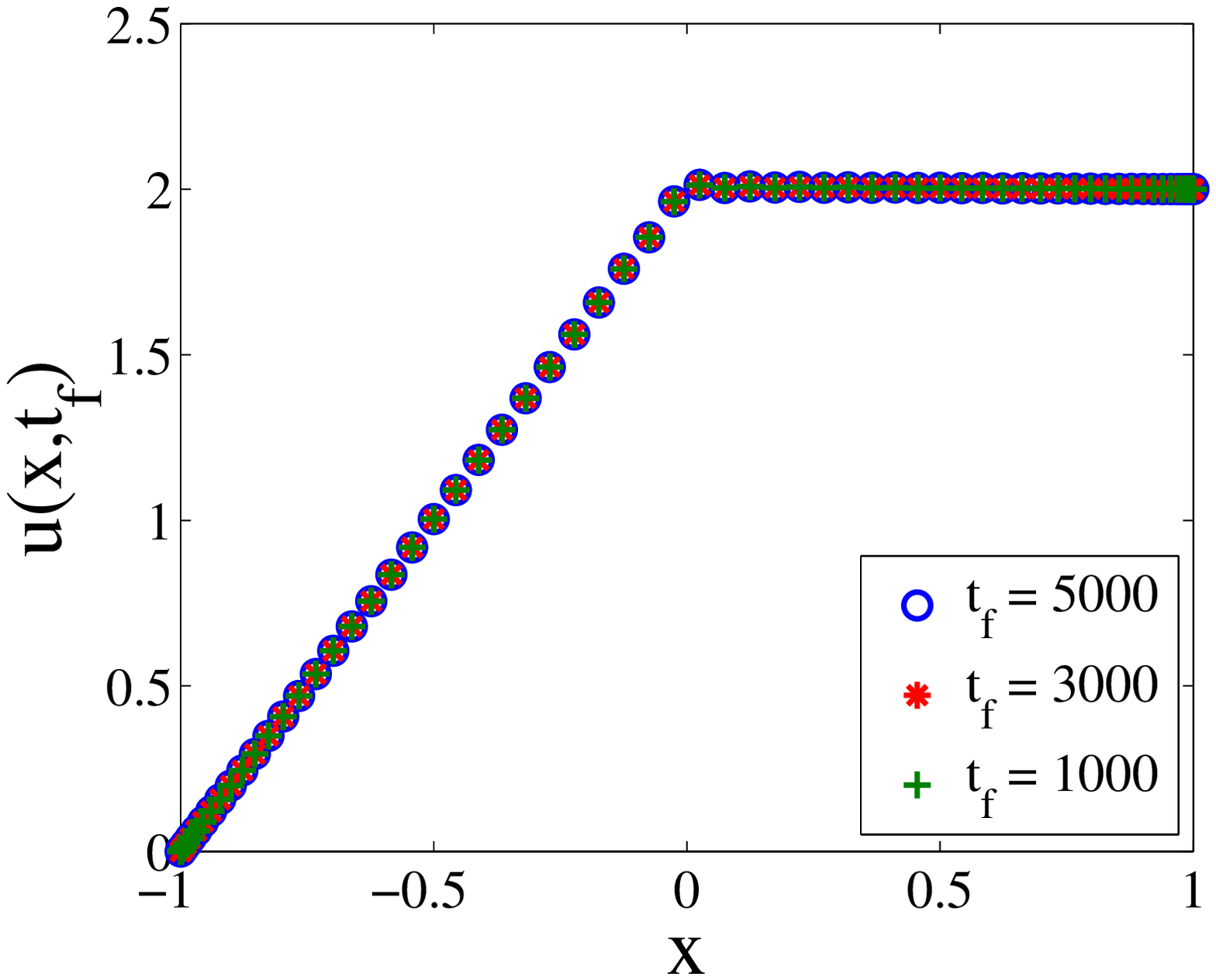}
		\includegraphics[scale=0.2]{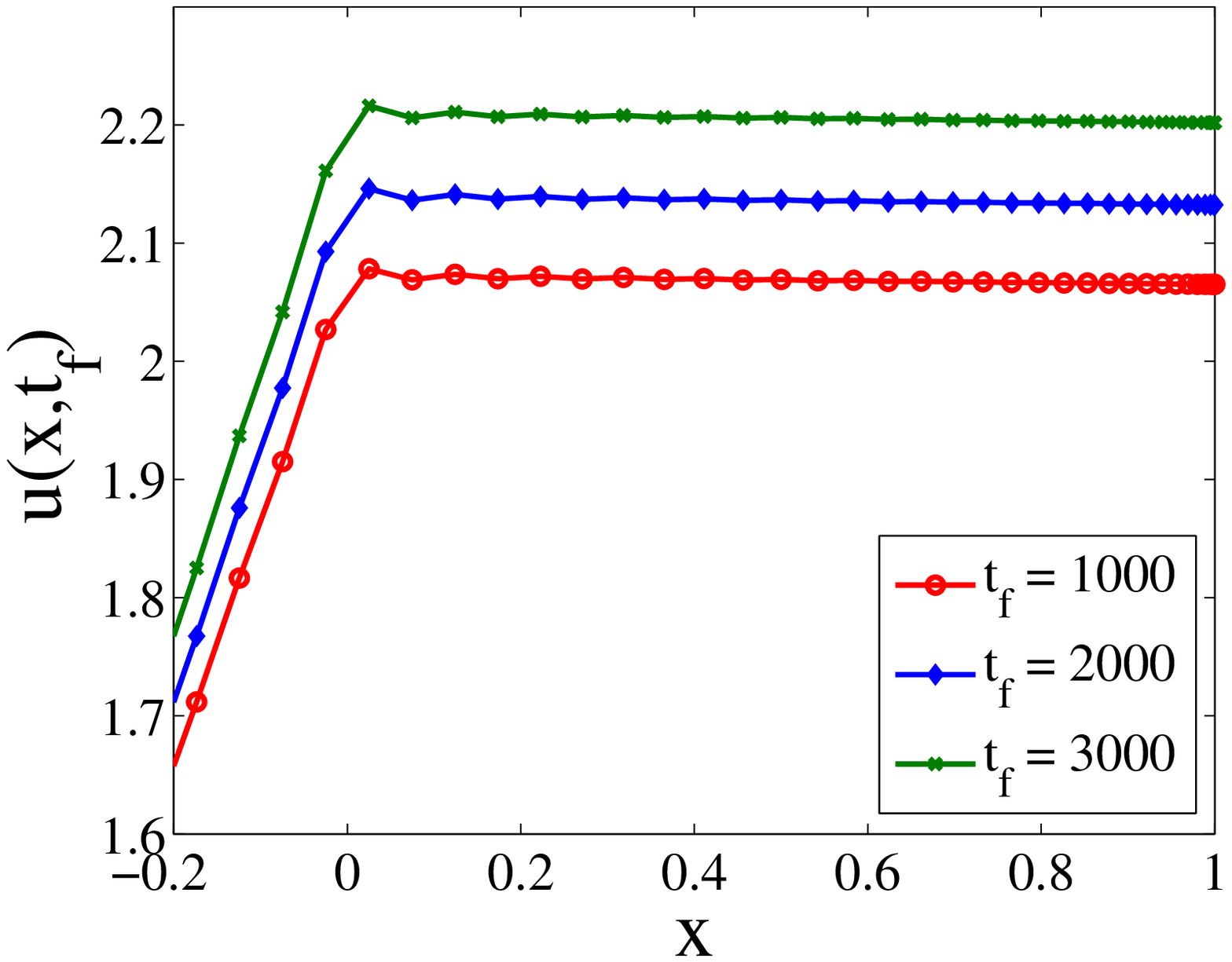}
		\includegraphics[scale=0.2]{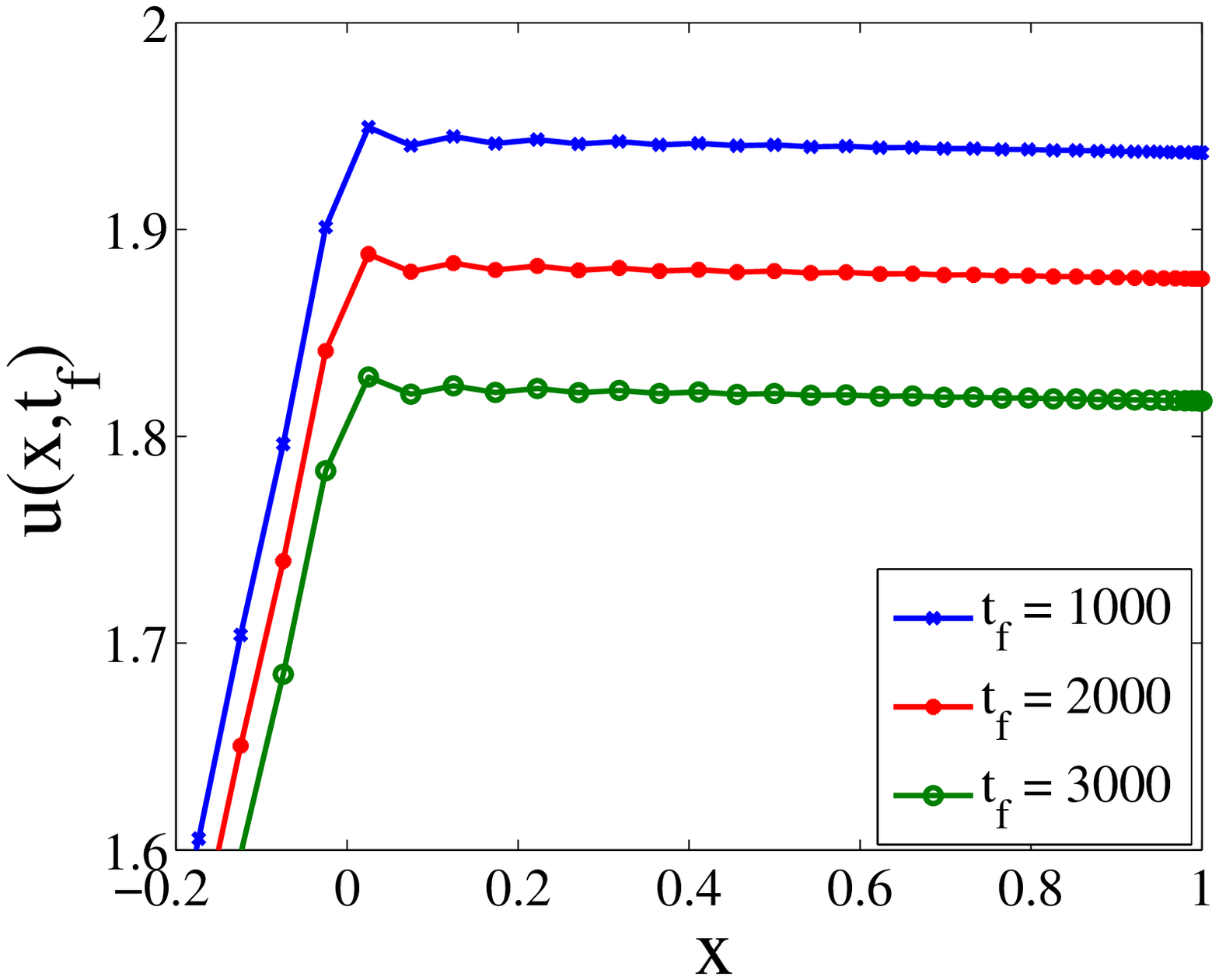}
		\caption{$u(x,t)$ for $\eta = 1.004332$ (left),  $\eta = 1.004364$ (middle),  $\eta = 1.004330$ (right). }
	 \label{fig:steady_states}
\end{figure}

In Figure \ref{fig:Eta} , we show  $\dot{E}(t)$ versus time with $\eta = \left\lbrace 1.004332,\; 1.004330,\; 1.004334 \right\rbrace $. We find for $\eta = 1.004332$ (blue), $|\dot{E}(t)| \lesssim O(10^{-6})$ when $t > 3000$. This shows $E(t)$ is almost constant with time. Another view of the steady-state solution is given in the right figure of Figure \ref{fig:Eta}, where some wiggles are shown due to the inconsistency of the initial condition in  the beginning but later the steady-state solution is obtained.   Note that we use the initial condition $u(x,0,\xi) = 1+x$. To verify the claim, we also calculate $u(x,t)$ with $\eta = 1.004332,\; 1.004300, \; 1.004364$ with different final time $t_{f}$ in  Figure \ref{fig:steady_states}. For $\eta = 1.00432, \; t_{f} = 1000, \; 3000,\; 5000,$ for $\eta = 1.004364\; \mathrm{and} \; \eta = 1.004330, \; t_{f} = 1000,\; 2000, \; 3000.$  The figure clearly shows that the solution remains almost constant for $\eta = 1.004332$. But for $\eta = 1.004300$ the solution decreases monotonically and for $\eta = 1.004364$ the solution increases monotonically.

%

\subsection{Finding the mean value}
First we compute the mean value of the solutions with time. Let ${\bar u}(x,t)$ be the mean over $\xi$. Then 
\begin{equation}
   {\bar u}(x,t) = {1\over 2}\int^1_{-1} u(x,t,\xi) d\xi. 
   \label{mean_int}
\end{equation}
We consider three mean values obtained by the following three methods. 

\textbf{{\em MC mean:}} To compute the above integral we use the uniform sampling of $\eta$ such that 
$$
            \eta^- = \eta_0 < \eta_1 < \cdots< \eta_i < \cdots < \eta_M = \eta^+
$$
and use the transformation $\xi = \frac{2\eta}{\eta^{+} - \eta^{-}}-\frac{\eta^{+} + \eta^{-}}{\eta^{+}-\eta^{-}}.$
Then  Eq.(\ref{mean_int}) takes the form 
$${\bar u}(x,t) = {\frac{1}{\eta^{+} - \eta^{-}}}\int^{\eta^{+}}_{\eta^{-}} u(x,t,\xi(\eta)) d\eta.$$
We use the trapezoidal rule to calculate the integration. \\

\textbf{{\em gPC mean:}}

From Eq.(\ref{mean_int}) we have 
\begin{eqnarray}
\bar{u}(x,t) = \frac{1}{2}\int_{-1}^{1}\sum_{k=0}^{N} \hat{u}_{k}(x,t)L_{k}(\xi)d\xi 
            = \frac{1}{2}\frac{2}{2l+1}\hat{u}_{l}(x,t)\delta_{0l}=\hat{u}_{0}(x,t). \nonumber
\end{eqnarray}
which is the first PC mode.\\

\textbf{{\em Quadrature mean:}} We consider $(N+1)$ Gauss-Legendre (GL) quadrature points $\eta_{i}$ between $[\eta^{-},\;\eta^{+}]$ and weights $\omega_{i}$. To get mean value we have 
$$\bar{u}(x,t) = \sum_{i=0}^{N}u_{i}(x,t_{f},\eta_{i})\omega_{i}.$$

\begin{figure}[h]
	\centering
		\includegraphics[width=0.8\textwidth]{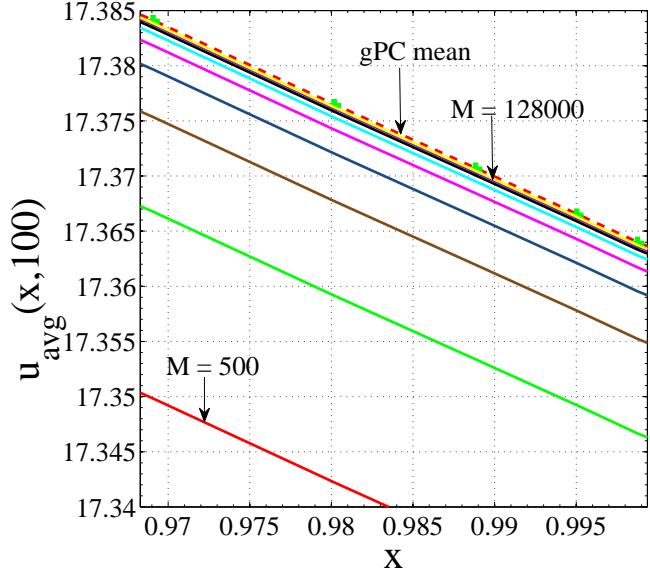}
		\caption{Convergence of $\bar{u}(x,t)$ by the MC simulations to gPC mean with different sample points.  The top curve (red dotted line) is the gPC mean. For MC simulations the number of sample points of $\eta$ are $M = 500\times 2^{j}, \;\mathrm{where}\; j = 0, \cdots ,8 $ .}
	 \label{fig:conv_mc}
\end{figure}

We calculate the mean value of the solution $\bar{u}(x,t,\eta)$ by the MC method at the final time $t_{f} = 100$ with different  uniform sample points $500$, $1000$, $2000$, $4000$, $8000$, $16000$, $32000$, $64000$ for $\eta \in [0.95, 1.05]$. Also we calculate the gPC mean at $t_{f} = 100$. The results are shown in Figure \ref{fig:conv_mc}. The MC mean solution converges to the gPC mean solution but the convergence is slow.  As shown in the figure, more than $128000$ MC samples are needed to match the gPC mean.  
This is because the energy growth/decay rate around the critical value of $\eta^c$ is extremely asymmetric and the uniform sampling is obviously inefficient. 

\begin{figure}
	\centering
		\includegraphics[width=0.49\textwidth]{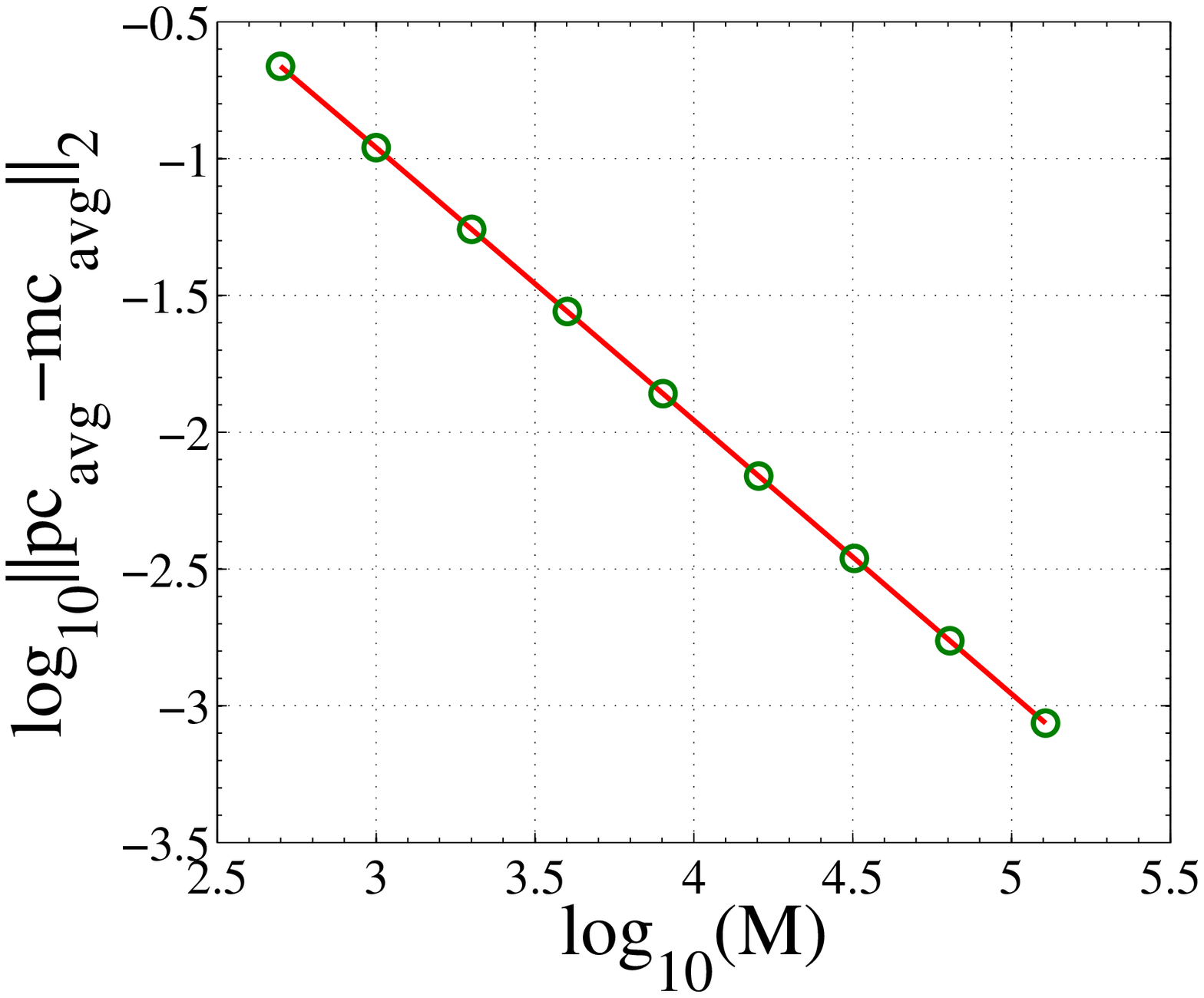}
		\includegraphics[width=0.45\textwidth]{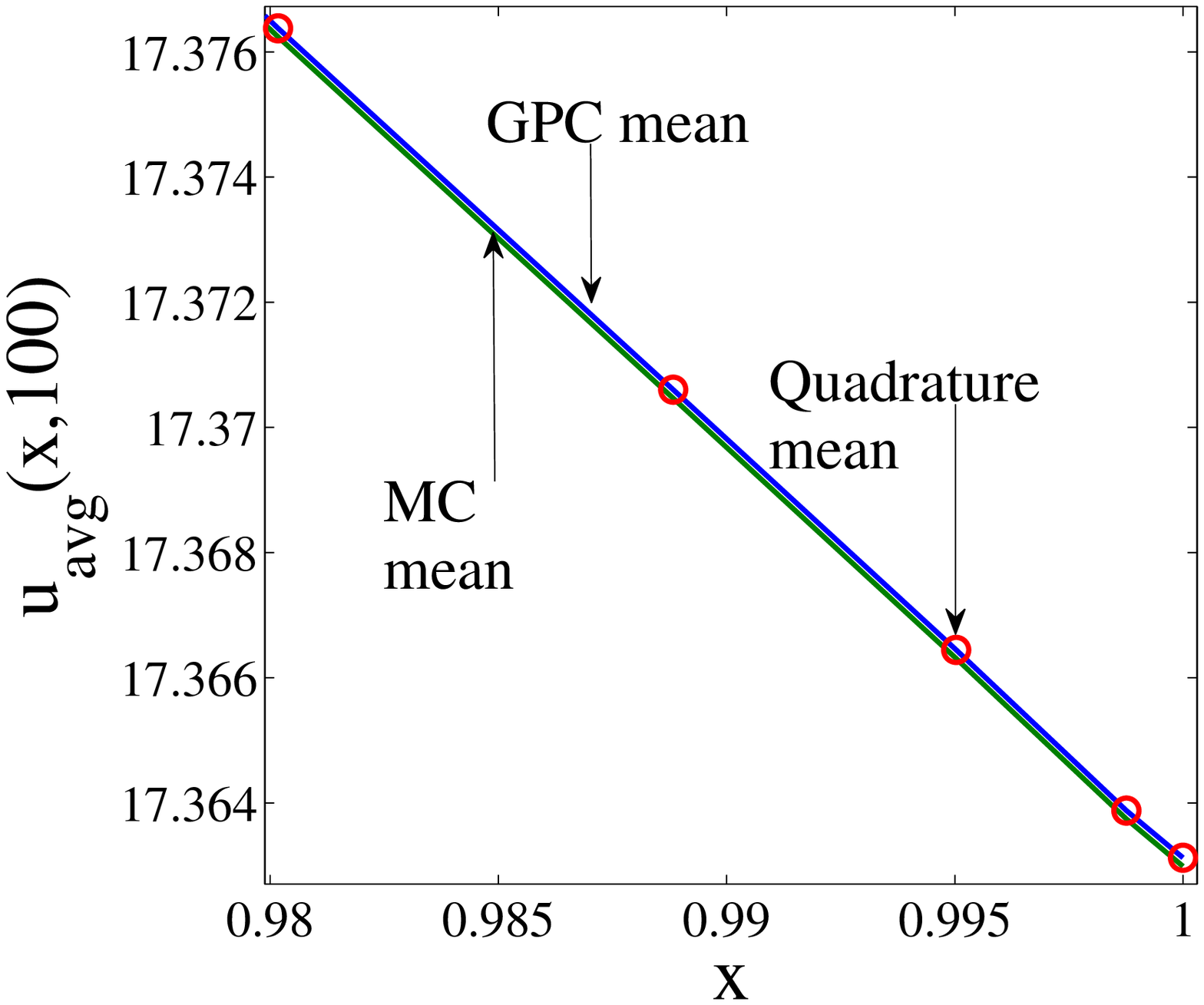}
		\caption{Left: Errors of the MC mean with different number of sample points $M$, where $M= 500 \times 2^{j}$ for $j = 0, 1, 2, \cdots , 8$.  
		Right: Comparison of $\bar{u}(x,100)$ derived from the MC simulation with $128000$ sample points, the Gauss-Legendre quadrature  with $50$ quadrature points (red circles) and the Legendre chaos method. It is observed that the quadrature mean and the PC mean are almost coincident.
		}
	 \label{fig:mcerror}
\end{figure}

The left figure of Figure \ref{fig:mcerror} shows the convergence of the MC simulations. We calculated the mean solutions at $t_{f}=100$ with the sample points $M_{j} = 500 \times 2^{j}$ where $j = 0, \cdots, 8$. We plot $||E_{j}||_{2}$ versus $M_{j}$  and  $E_{j} = PC_{avg}-\bar{u}_{j}$. Here $\bar{u}_{j}$ is the mean solution of  uniformly distributed $M_{j}$ samples. The result is  plotted in logarithmic scale and  the decay is linear.

The right figure of Figure \ref{fig:mcerror} shows the mean solutions with the Legendre gPC, MC simulation with $128000$
equally distributed sample points  and the quadrature mean with $50$ Gauss-Legendre quadrature points in $[0.9, 1]$. The solid red line represents the result by the  MC simulation which deviates from the gPC mean or GL quadrature mean. The GL quadrature mean with the small number of points coincides with PC mean. Using the quadrature sampling for the mean solution was similarly used in \cite{Dalbey}.

\subsection{Finding critical value of $\eta$} 

We first define a solution sequence $U^k(\eta)$ at $x = 1$. The index $k$ denotes the sequence corresponding to time $t_k$ such as 
$$
               k: t_k = k \Delta t, \quad \Delta t > 0.
$$
Define $\eta^l$, $l = 1, 2, 3, \cdots$, 
$$
    \eta^l = \left\{ \eta| U^l(\eta) - U^{l-1}(\eta) = 0\right\}. 
$$
For large $t$, $U^k(\eta)$ is a monotone function with respect to $\eta$ around the critical value of $\eta$. As there exists the critical value of $\eta$, $\eta_c$ , there exists $\eta$ such that $U^l(\eta) - U^{l-1}(\eta) = 0$.
Now we define another sequence $\Delta \eta$ such as
$$
    \Delta \eta^l = \eta^l - \eta^{l-1}. 
$$ 
As $U^l$ is monotone with respect to $\eta$ for large $t$, we have
$$
   \lim_{l \rightarrow \infty} \Delta \eta^l = 0.  
$$
and 
$$
  \lim_{l \rightarrow \infty} \eta^l = \eta_c. 
$$

\begin{figure}
	\centering
		\includegraphics[scale = 0.3]{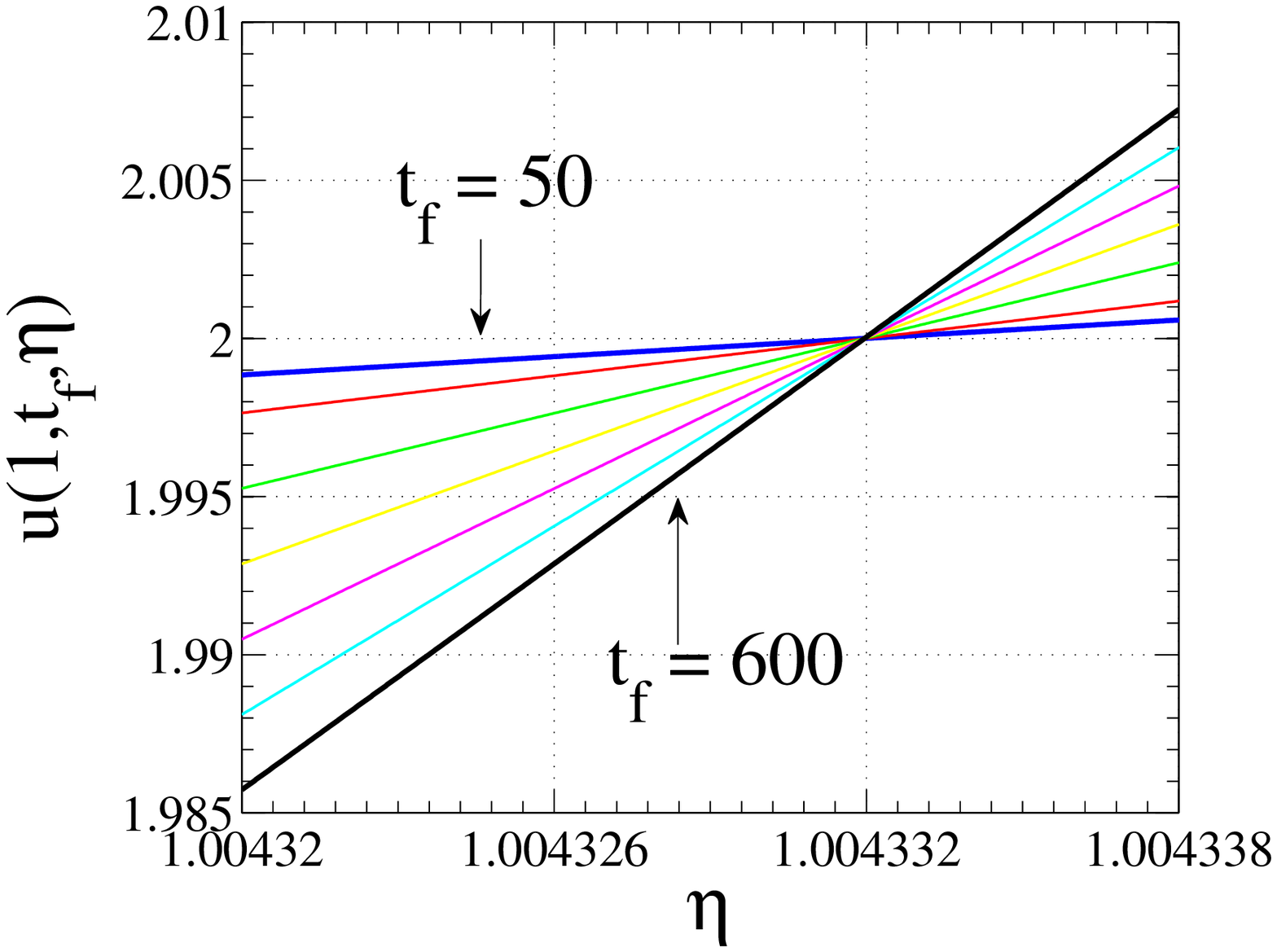}    
		\includegraphics[scale = 0.3]{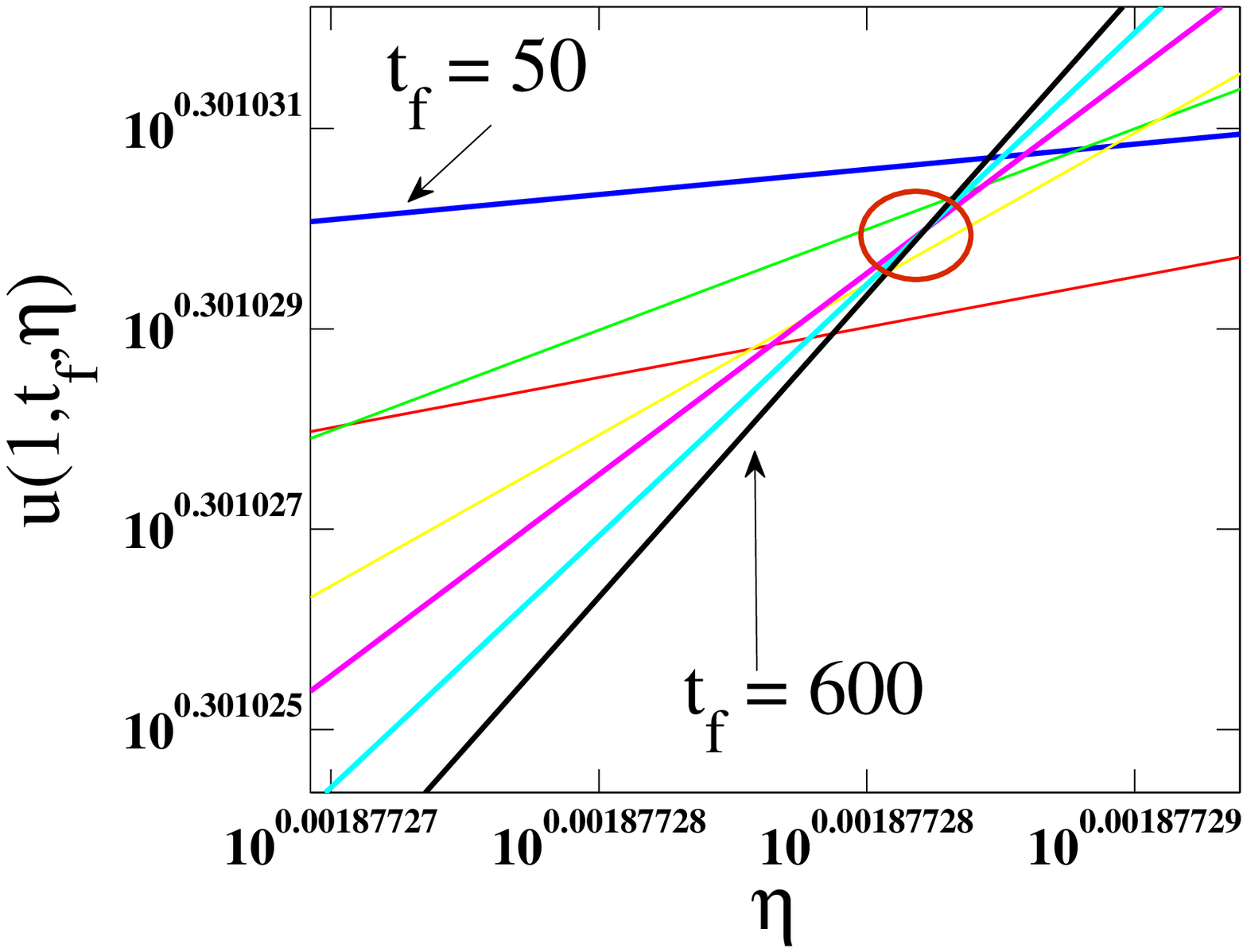}
		\caption{ Top: The locus of $u(1,t_{f},\eta)$ with different final time by the MC simulations. $\eta$ varies from $0.95$ to $1.05$. The curves seem to intersect at the common point.  The lowest curve (blue solid line) in the right side of intersection point is for $t_{f} = 50$ and then the top curve (black solid) for $t_{f} = 600$. 
		Bottom: The close view  around the intersection point. For the MC simulations, the intersection points are distributed randomly in a very small region but for higher time the curves tend to converge to a single point.}
	 \label{fig:crit_eta_mc}
\end{figure}

Figure \ref{fig:crit_eta_mc} shows the sequence of $U^{l}$ and $\eta^{l}$ by the MC simulations with $10^5$ sample points at different time intervals. The final time varies from $50$ to $600$ and we calculate the solution at $t_{f} = 50$, $100$, $200$, $\cdots$, $600$. For any fixed $t_{f}$, we plot the  solution for different $\eta$. If $\exists \: \eta_{c}$   then the graph $u(1,t_{f},\eta)$ converges to zero as 
$t \rightarrow \infty$ when $\eta <  \eta_{c}$ and similarly $u(1,t_{f},\eta)$ diverges to $\infty$ as $t \rightarrow \infty$ when $\eta > \eta_{c}$. From Figure \ref{fig:crit_eta_mc}, we find that the solutions intersect at the common point $(\eta^{n}_{c}, u_{s}(1,t_{f},\eta^{n}_{c}))$ and we claim that such $\eta^{n}_{c}$ is the critical value of $\eta$ and $u_{s}(1,t_{f},\eta^{n}_{c})$ is  the steady-state solution.  

The bottom  figure  shows  the magnified region near $\eta^{n}_{c}$. 
If we see the figure carefully, it is then observed that the solutions for $t_{f} = 400$,  $ t_{f} = 500$ and  $t_{f} = 600$ intersect almost at the same point, whose coordinate is approximately equal to $( 10^{0.00187728},10^{0.3010295}) = (1.004331952841278,   1.999997717384314)$. Then our estimates are $\eta^{n}_{c}  = 1.004331952841278 $ and $u_{s}(1, t_{f}, \eta^{c}_{n}) = 1.999997717384314$.
In Section 4 we found for $(m+1) = 64$, the critical value of $\eta$ is about $1.004332$, which is close
to the value of $\eta^n_c$. Here note that the value of $\eta^c$ in Table 2 was obtained by many of the individual simulations of the equation with different value of $\eta$, which is highly expensive and $\eta^c$ was obtained with only $10^-6$ accuracy. 
The error for the steady-state solution at $x = 1$ is $2 - 10^{0.3010295} \approx 2.2 \times 10^{-6}$, which agrees our error estimation in Section 5.1 and Figure \ref{fig:Eta}.  The error of $10^{-6}$ comes from the Gibbs oscillations at $x = 1$ as explained. 

\begin{figure}
	\centering
		\includegraphics[scale = 0.3]{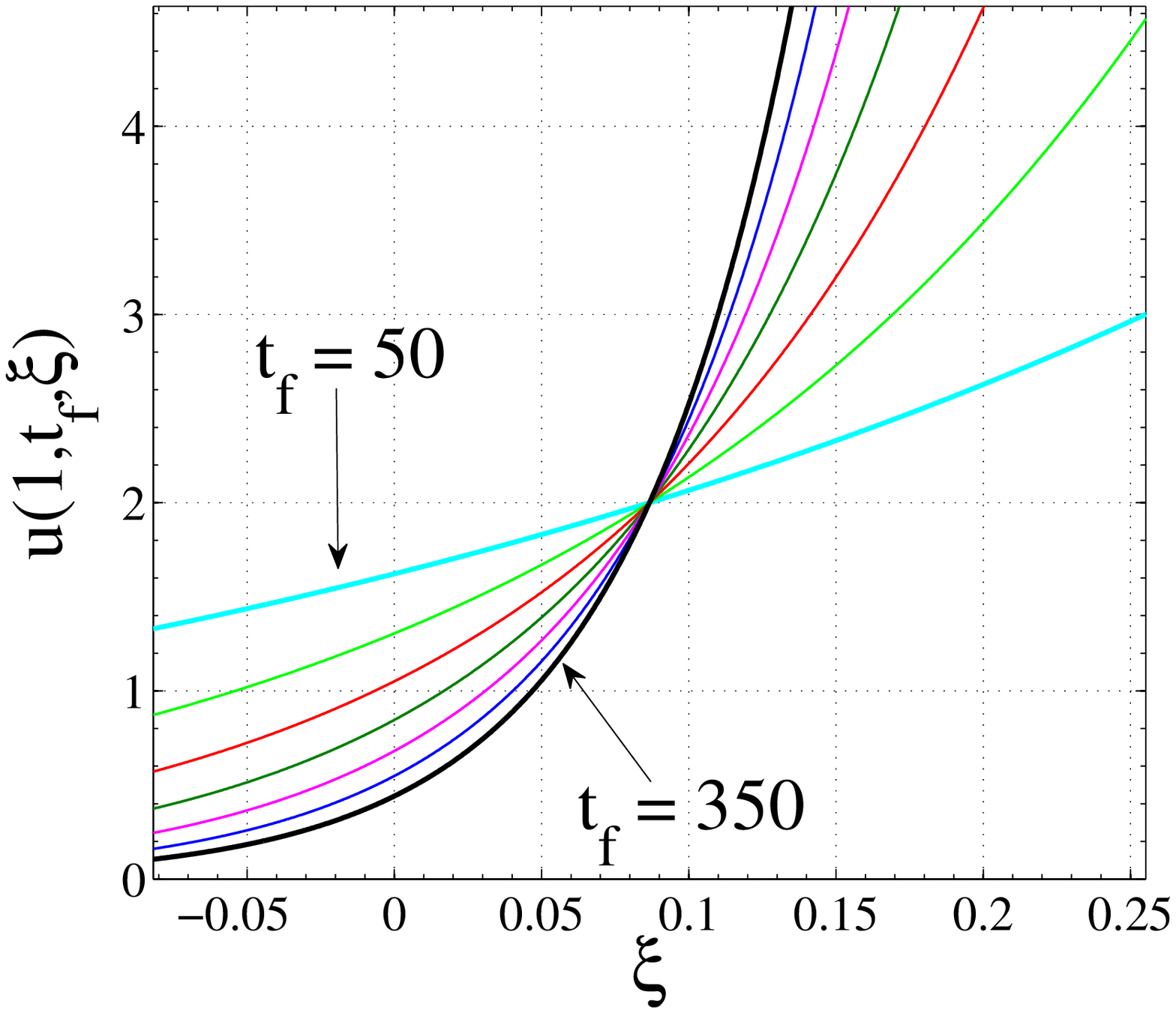}   
		\includegraphics[scale = 0.3]{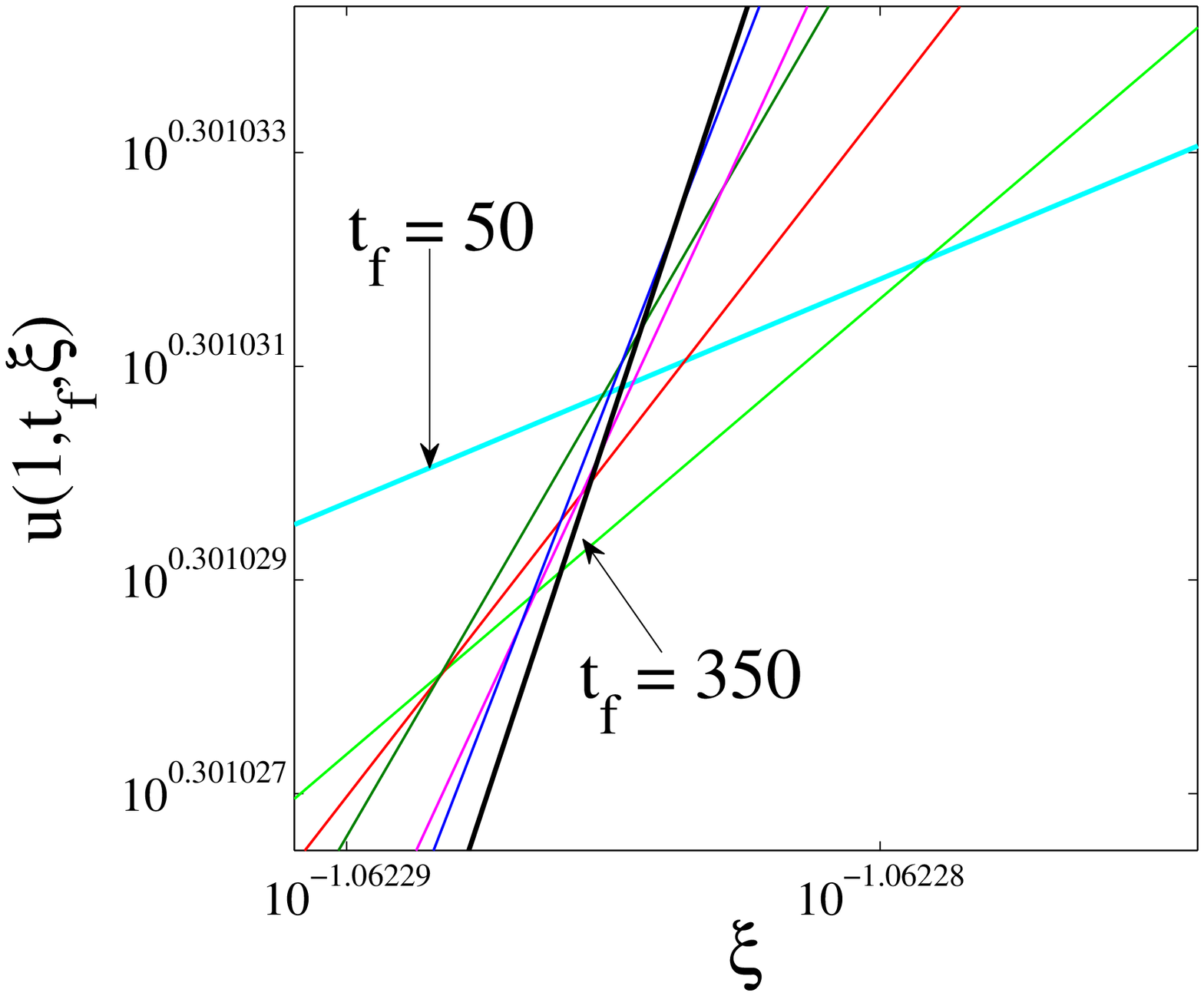}
		\caption{ Left: The locus of $u(1,t_{f},\eta)$ with different final time by the Legendre chaos. $\eta$ varies from $0.95$ to $1.05$.The curves intersect at the common point.  The lowest curve (cyan solid line) in the right side of intersection point is for $t_{f} = 50$ and the top curve (black solid) for $t_{f} = 500$. 
		Right: A close view  around the intersection point. For the Legendre chaos method, the intersection points are distributed randomly in a very small region. Note that the horizontal axis is for $\xi$ not $\eta$.}
	 \label{fig:crit_eta_pc}
\end{figure}


Figures \ref{fig:crit_eta_pc}  shows how to find the critical value $\eta_{c}$ by the gPC method. In the gPC method we reconstruct the solutions by considering the first $80$ modes, at different time $t_{f}$, where $t_{f}$ runs from $50$ to $350$ with the interval length $50$. At every stage we plot $u(1,t_{f}, \xi)$, where $\xi \in [-1, 1]$ and consider $10^{6}$ equally distributed sample points generated by the $80$ gPC modes. We found that the solutions intersect at the common point $(\xi_{c}, u_{s}(1,t_{f},\xi_{c}))$ and we estimate that it is the  critical value of $\xi$ for $\eta \in \left[0.95,\; 1.05 \right]$. 
We find the approximate value of $\xi_{c}$ is  $\xi_{c}= 10^{-1.062285} \approx 8.6639313 \times 10^{-2}$ and $u_{s}(1, t_{f},\xi_{c})=10^{0.301032} \approx 2.000009230329$. For this value of $\xi_{c}$ the corresponding $\eta_{c}$ is $\eta_{c} = 1.004331965650747$. The value of $\eta_{c}$ divides the interval $[0.95,\; 1.05]$ into two subintervals $[0.95, \eta_{c}]$ and $[\eta_{c}, 1.05]$. When $\eta < \eta_{c}$, the graph converges to zero as $t_{f} \rightarrow \infty$ and similarly the graph diverges to $\infty$ if $\eta > \eta_{c}$ as $t_{f} \rightarrow \infty$.
In Figure \ref{fig:crit_eta_pc}, the right figure is the magnified version of the left figure near $\xi_{c}$. 
The graphs intersect within a very small region between $10^{0.301028} \approx 1.999990809648848$ and $10^{0.301032} \approx 2.000009230329776$ in the vertical axis and between $10^{-1.06229} \approx 8.663831554873760 \times 10^{-2}$ and $10^{-1.06228} \approx 8.664031049264383 \times 10^{-2}$ in the horizontal axis.  The last two lines at $t_{f} = 300$ and $t_{f} = 350$ (blue and black) intersect at $u_{s} = 10^{0.301032} \approx 2.000009230329776$ and the error in this case $\approx 7 \times 10^{-6}$. The error in calculating the critical value of $\eta$ must be less than the difference of $10^{-1.06229}$ and $10^{-1.06228}$ which is almost equal to $1.99 \times 10^{-6}$. Thus we get the result with in the accuracy of $\sim 10^{-6}$. The difference between $\eta^c$ computed by the MC method and the gPC method is
$$
       \Delta \eta^c :=| \eta^c (gPC) - \eta^c (MC) | \sim  1.2809 \times 10^{-8}. 
$$
Here note that for the gPC method we only use $N = 80$ while for the MC method we use $M = {10}^5$, which shows that the gPC method can determine $\eta^c$ efficiently. 

\section{Critical phenomenon of the sine-Gordon equation}
\subsection{Basic idea and problem formulation}
Consider the sine-Gordon equation 
$$u_{xx} - u_{tt}  + \sin(u) = 0, \quad x \in \mathbb{R},$$
where $u:\mathbb{R}\times\mathbb{R}^+ \rightarrow \mathbb{R}$. 
The \textit{soliton (or kink)} solution of the sine-Gordon equation is given by 
$$ u(x,t) = 4\tan^{-1}\left\lbrace \exp\left(\frac{x-Vt}{\sqrt{1-V^2}} \right)  \right\rbrace. $$
This represents the continuous profile with $u \rightarrow 0$ as $x \rightarrow -\infty$, and $u \rightarrow 2\pi$ as $x \rightarrow +\infty$. This solution propagates with the constant velocity $V$. Another solution can be obtained in the form 
$$ u(x,t) = 4\tan^{-1}\left\lbrace \exp\left(-\frac{x-Vt}{\sqrt{1-V^2}} \right)  \right\rbrace = 4\cot^{-1}\left\lbrace \exp\left(\frac{x-Vt}{\sqrt{1-V^2}} \right)  \right\rbrace, $$
which is known as  the \textit{antisoliton(or antikink)}.

We consider the sine-Gordon equation with an impurity
\begin{eqnarray}
\label{SG1}
u_{tt} - u_{xx} + \sin(u) = \epsilon\delta(x)\sin(u), \quad x \in \mathbb{R} ,\quad t \in \mathbb{R}^{+},
\end{eqnarray}
where $\epsilon \in \mathbb{R}^+$ is the amplitude or strength of the singular potential function. 

The full sine-Gordon equation with an impurity Eq. (\ref{SG1}) is, however, non-integrable and no exact solution is known. The problem we solve numerically in this work is defined by Eq. (\ref{SG1}) in $x \in \Omega$, where $\Omega = [-L, L]$, $L > 0$ with the initial kink velocity $V$. For simplicity we assume that the singular potential exists at $x = 0$. For the numerical computation, we set $L = 8$, $\epsilon = {1\over 2}$, and the initial kink front $x_0 = -6$. 
The initial and boundary conditions are given below. 
%
%
%
%

Initial conditions (kink soliton)
\begin{eqnarray}
u(x,0) &=& 4\tan^{-1}\left[ \exp \left\lbrace \frac{x-x_{0}}{\sqrt{1-V^2}} \right\rbrace \right], \nonumber \\
u_{t}(x,0) &=& -\frac{4V}{\sqrt{1-V^2}}\left[\frac{\exp\left(\frac{x-x_{0}}{\sqrt{1-V^2}} \right) }{1+\exp\left(\frac{2(x-x_{0}}{\sqrt{1-V^2}} \right) } \right].  \nonumber
\end{eqnarray}

Boundary conditions:
\begin{eqnarray}
u_{t}(L,t) +u_{x}(L,x) &=& 0,  \nonumber \\
u(-L,t) &=& 4\tan^{-1}\left\lbrace \exp\left[\frac{-L-x_{0}-Vt}{\sqrt{1-V^2}} \right] \right\rbrace . \nonumber
\end{eqnarray}

The above boundary condition is used for the outflow boundary condition. Note that this boundary condition is not an exact outflow boundary condition because of the nonlinear potential term, $\sin(u)$. Better boundary conditions and exact outflow boundary conditions can be found in \cite{Zheng}. Since the main goal of our study is to find the critical velocity, the above boundary condition does not affect the results.  

\begin{figure}
  \centering
    \includegraphics[width=0.4\textwidth]{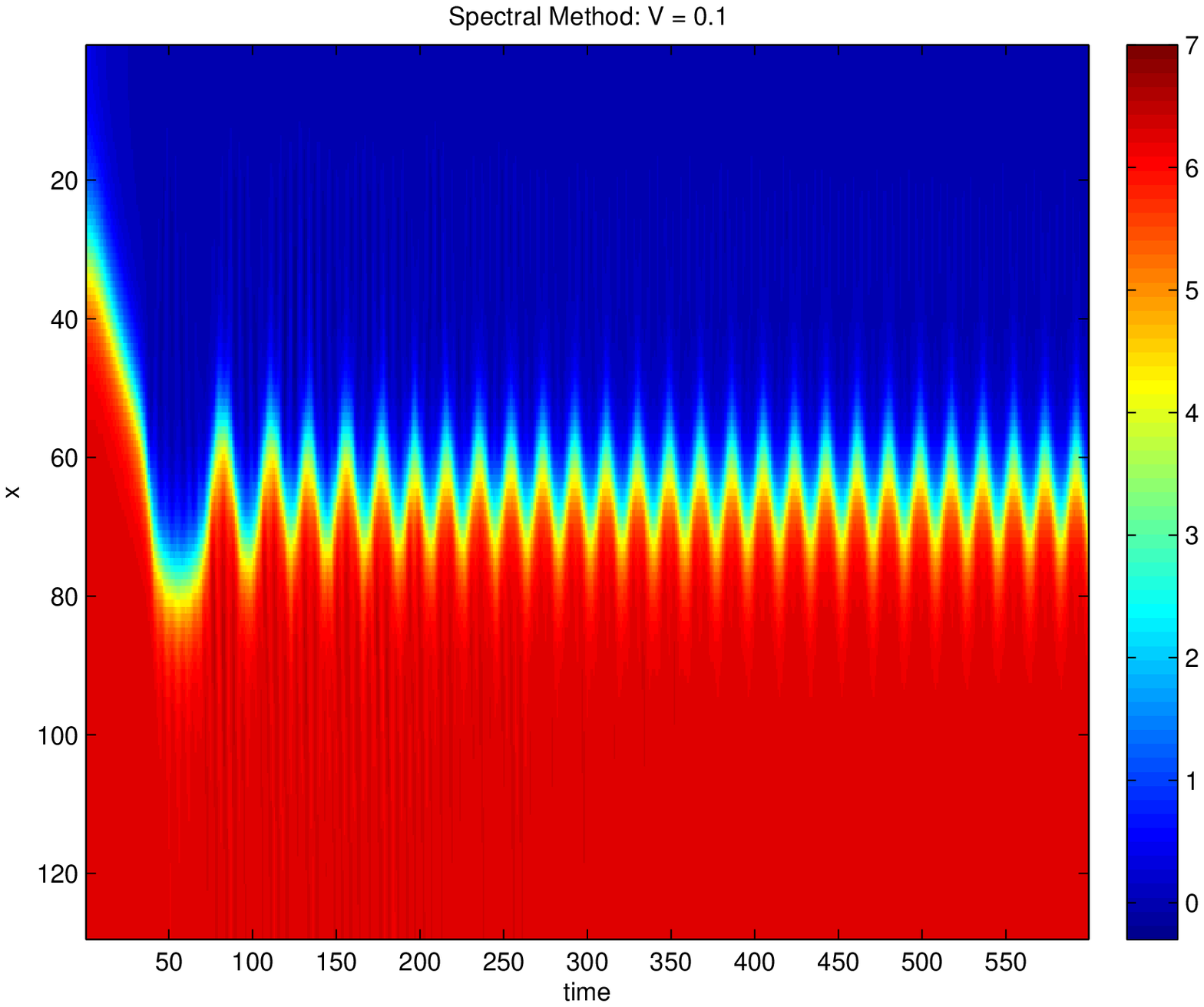}
  \includegraphics[width=0.4\textwidth]{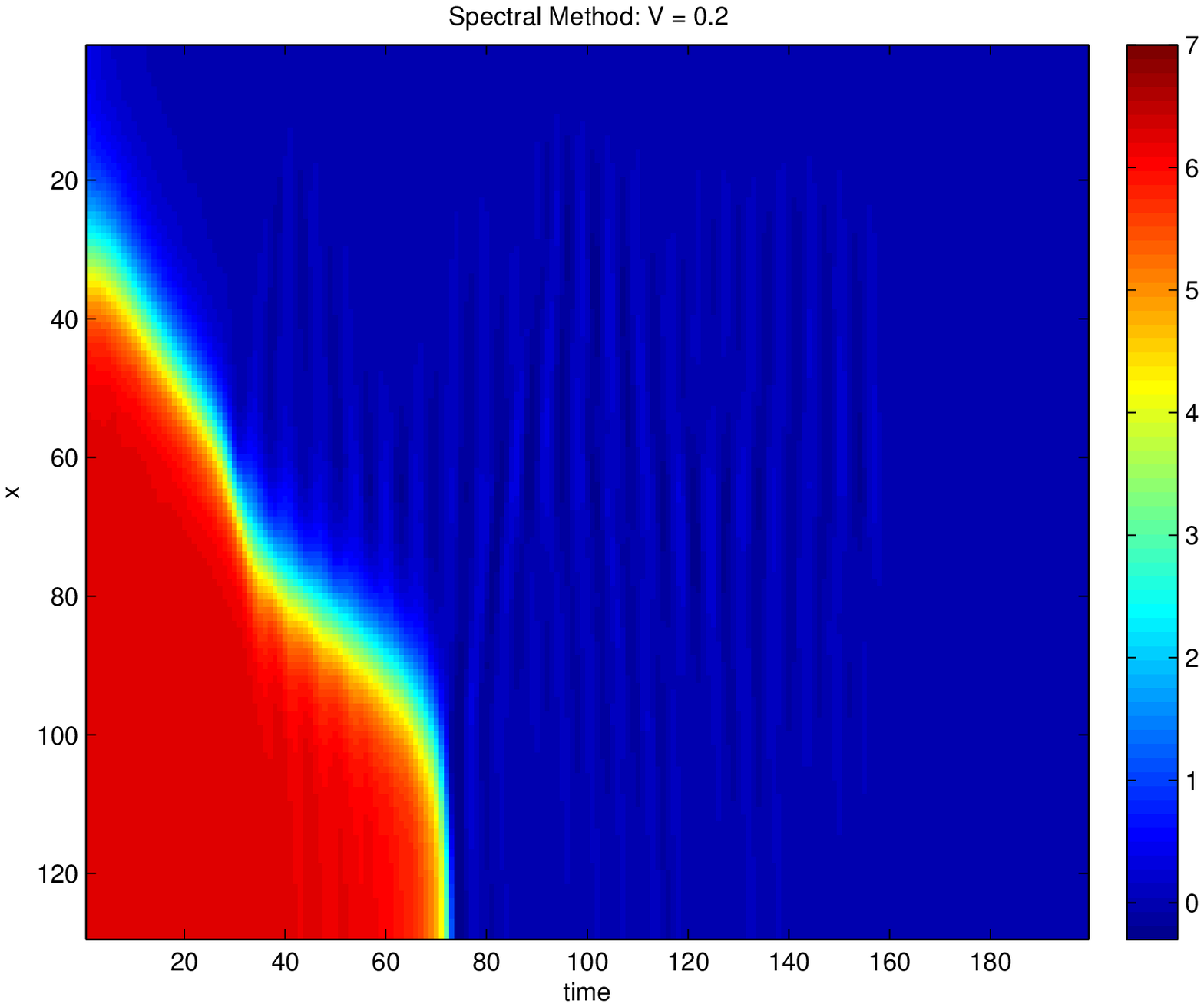}
  \caption{Spectral solution of the sine-Gordon equation with $m = 127$. Left: Sub-critical behavior with $V = 0.1$. Right: Super-critical behavior $V = 0.2$ \cite{JungDon2009}.}
\label{critical}
\end{figure}

Figure \ref{critical} shows the solution behaviors around the critical value of $V_c$, $0.1 < V_c < 0.2$. The left figure shows the sub-critical solution behavior when $V = 0.1< V_c$ and the right figure when $V = 0.2 > V_c$. The horizontal axis denotes the time and the vertical axis the location of the kink front (the red color shows the kink front and the $2\pi$ plateau). These solutions are obtained by the direct numerical approximation using the Chebyshev collocation method with the consistent approximation of the $\delta$-function explained in the previous section with $m+1 = 128$. As shown in the figure, the soliton solution is captured in the potential well if $V< V_c$ while it eventually passes through the potential well. That is, at $x = L$, $u(x,t)$ becomes either $2\pi$ (particle-capture or particle-reflection) or $0$ (particle-pass) as $t \rightarrow \infty$. We use this solution behavior to find the critical velocity $V_c$. Here note that the early numerical study of the discretized sine-Gordon equation is found in the literature \cite{PK1984}.  

\subsection{Finding the critical velocity}


In the previous section, we considered the Klein-Gordon type linear equation for which we use the MC sampling using the gPC expansion for constructing the solutions. However, for the nonlinear sine-Gordon equation, such a MC sampling is not usable because of the two extreme PDF of the solution for large $t$.  The construction of the gPC for large $t$ yields either \textit{trapped} or \textit{transmitted} solution. 

\begin{small}
\setlength{\unitlength}{4cm}
\begin{picture}(1,1)

\put(0.5,0){\vector(1,0){1.3}}
\put(0.3,0.3){$f(z)$}
\put(0.5,0){\vector(0,1){0.7}}
\put(0.7,0){\line(0,1){0.6}}
\put(1.4,0){\line(0,1){0.4}}
\put(0.7,-0.1){$0$}
\put(1.4,-0.1){$2\pi$}
\put(1.0,-0.2){$z = u^\infty$}
\put(1.0,0.53){\vector(-1,0){0.25}}
\put(1.05,0.5){$(1-V_c)\delta(t)$}
\put(1.7,0.23){\vector(-1,0){0.25}}
\put(1.75,0.2){$V_c\delta(t-2\pi)$}
\put(0.05,-0.3){Diagram 1: The PDF of the solution for $t \rightarrow \infty,\; u^{\infty}$ at $x = L$ for uniform $V$.}
\end{picture}
\end{small}

%

\vskip .8in
The critical value can be computed by the weighted ratio of the two extreme PDF as shown in  Diagram 1 if we know the PDF of $V$. The ratio is indeed the amount of how much the mean value deviates from the value of $2\pi$ for the \textit{trapped} solution. Thus the critical value of $V$ can  be still found from the gPC expansion using the mean mode. 
%
%
%
If $V \in \left[V_{a},V_{b}\right]$ and if $V_{c}$ is the critical velocity, then for $V \in \left[ V_{a}, V_{c}\right]$, the solution is $2\pi$ (because the solution is \textit{trapped}), and for $V \in \left( V_{c}, V_{b}\right]$, the solution vanishes (because the solution is \textit{passed}). If $f\left(V; V_{a}, V_{b}\right)$ is the PDF for $V$, then
\begin{eqnarray}
\overline{u\left(x=L,t_{f}\right)} &=& \int_{V_{a}}^{V_{b}}f(V;V_{a},V_{b})u(x,t_{f},V)dV \nonumber\\
                 &=& 2\pi\int_{V_{a}}^{V_{c}}f(V;V_{a},V_{b})dV .
                 \label{critical_velo_pdf}
\end{eqnarray}
\noindent
For example, if we use the Legendre chaos, the distribution function is uniform and $f\left(V;V_{a},V_{b}\right)= \frac{1}{V_{b}-V_{a}}$. For this by using Eq. (\ref{critical_velo_pdf}), the Legendre mean $\bar{u}_{leg}$ is given by 
\begin{eqnarray}
\bar{u}_{leg} &=& \frac{2\pi}{V_{b}-V_{a}}\left(V_{c}-V_{a}\right), \nonumber
\end{eqnarray}
and we have 
\begin{eqnarray}
            V_{c}  &=& V_{a} +\left(V_{b}-V_{a}\right)\frac{\bar{u}_{leg}}{2\pi}.
            \label{critical_legendre}
\end{eqnarray}
%
%

\subsection{Hermite chaos and the distribution function of input velocity $V$}

Suppose that the input velocity, $V$, is the random variable which is normally distributed with mean $\mu$ and standard deviation $\sigma$. Here $V$, conditional on $\alpha \le V \le \beta$, has a truncated normal distribution. We consider $\mu = \frac{1}{2}\left(\alpha + \beta\right)$ and $\sigma \sim 10^{-2}$ where $\alpha = V_a, \beta = V_b$. The reason for choosing small $\sigma$ is to keep the error $\sim 10^{-6}$.  The PDF of $V$, $f$, is given by
\begin{equation}
f\left(V; \mu,\sigma,\alpha, \beta \right) = \left\{\begin{array}{cc} \frac{\frac{1}{\sigma}\psi\left(\frac{V-\mu}{\sigma}\right)}{\Phi\left(\frac{\beta - \mu}{\sigma}\right)-\Phi\left(\frac{\alpha - \mu}{\sigma}\right)} ,&\mathrm{for} \quad \alpha \le V \le \beta, \\
   0, & \mathrm{otherwise} . \end{array} \right.
   \label{hermite-pdf}
\end{equation}
\noindent
The density function involves $\psi(\cdot)$ which is the PDF of the standard normal distribution and $\Phi(\cdot)$, its cumulative distribution function \cite{Chen_Gottlieb}. Due to the small value of $\sigma$,
$$
\int_{\alpha}^{\beta}e^{-\frac{\left(V-\mu\right)^{2}}{2\sigma^{2}}}H_{l}\left(\frac{V-\mu}{\sigma}\right)H_{l^{'}}\left(\frac{V-\mu}{\sigma}\right)dV \approx \sigma\sqrt{2\pi}l!\;\delta_{ll^{'}},
$$
for which see Appendix 1. 
\subsection{Sine-Gordon equation and Hermite chaos}
We expand $u(x,t,V)$ in terms of the Hermite polynomials as follows
$$u(x,t,V) = \sum_{l=0}^{\infty} \hat{u}_{l}(x,t)H_{l}(\frac{V-\mu}{\sigma}).$$
For the numerical purpose, we consider the first $N+1$ modes.
Then the sine-Gordon equation (\ref{SG1}) can be written as 
\begin{eqnarray}
[\hat{u}_{l}(x,t)]_{tt} &=& [\hat{u}_{l}(x,t)]_{xx} + \frac{1}{l!\sqrt{2\pi}}\left\lbrace \epsilon \delta(x) -1 \right\rbrace \phi_{l}(x,t), \nonumber \\
 \mathrm{where}\qquad \phi_{l}(x,t) &=& \int_{-\infty}^{\infty}\sin\left(\sum_{l=0}^{N}\hat{u}_{l}(x,t)H_{l}(\gamma) \right)H_{l}(\gamma)\omega(\gamma)d\gamma, \;  \nonumber
\end{eqnarray}
and
$$
\gamma = \frac{V-\mu}{\sigma}.
$$
For the initial condition, we have
$$4\tan^{-1}\left[ \exp \left\lbrace \frac{x-x_{0}}{\sqrt{1-V^2}} \right\rbrace \right] =\sum_{l=0}^{\infty} \hat{u}_{l}(x,0)H_{l}\left(\frac{V-\mu}{\sigma}\right). $$
By using the orthogonal property of the Hermite polynomials we get
$$\hat{u}_{l}(x,0) = \frac{4}{\sqrt{2\pi}l!}\int_{-\infty}^{\infty}\tan^{-1}\left[ \exp \left\lbrace \frac{x-x_{0}}{\sqrt{1-\left(\mu+\sigma\gamma\right)^2}} \right\rbrace \right]H_{l}(\gamma)\omega(\gamma)d\gamma.$$
Similarly 
$$\left(\hat{u}_{l}(x,0) \right)_{t} = \frac{-4}{\sqrt{2\pi}l!}\int_{-\infty}^{\infty}\frac{\left(\mu+\sigma\gamma\right)}{\sqrt{1-\left(\mu+\sigma\gamma\right)^2}}\frac{\exp(\frac{x-x_{0}}{\sqrt{1-\left(\mu+\sigma\gamma\right)^2}})}{1+\exp\left(\frac{2(x-x_{0})}{\sqrt{1-\left(\mu+\sigma\gamma\right)^2}}\right)}H_{l}(\gamma)\omega(\gamma)d\gamma.$$

\begin{figure}
	\centering
		\includegraphics[height=0.43\textwidth]{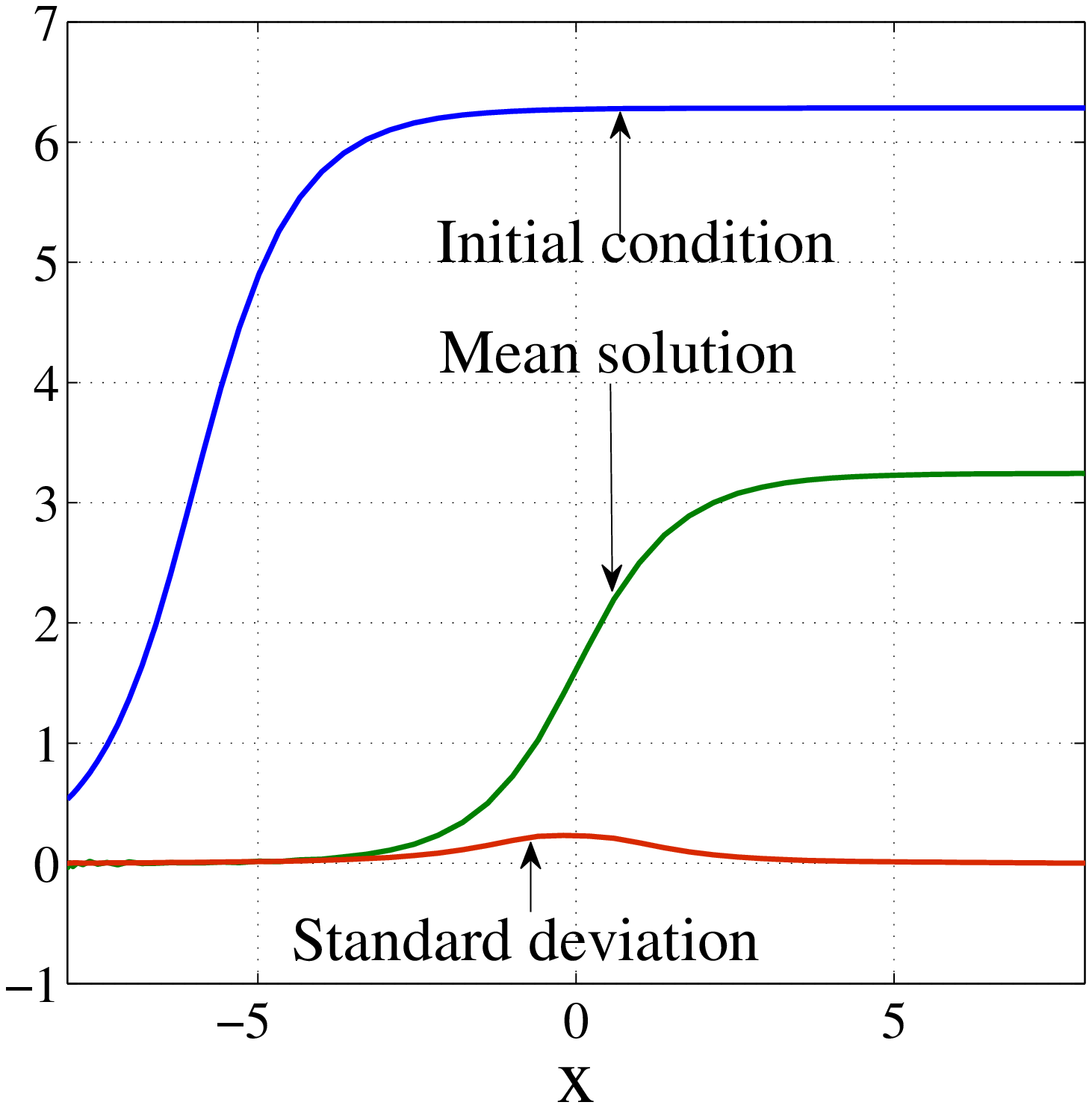}    
		\includegraphics[height=0.43\textwidth]{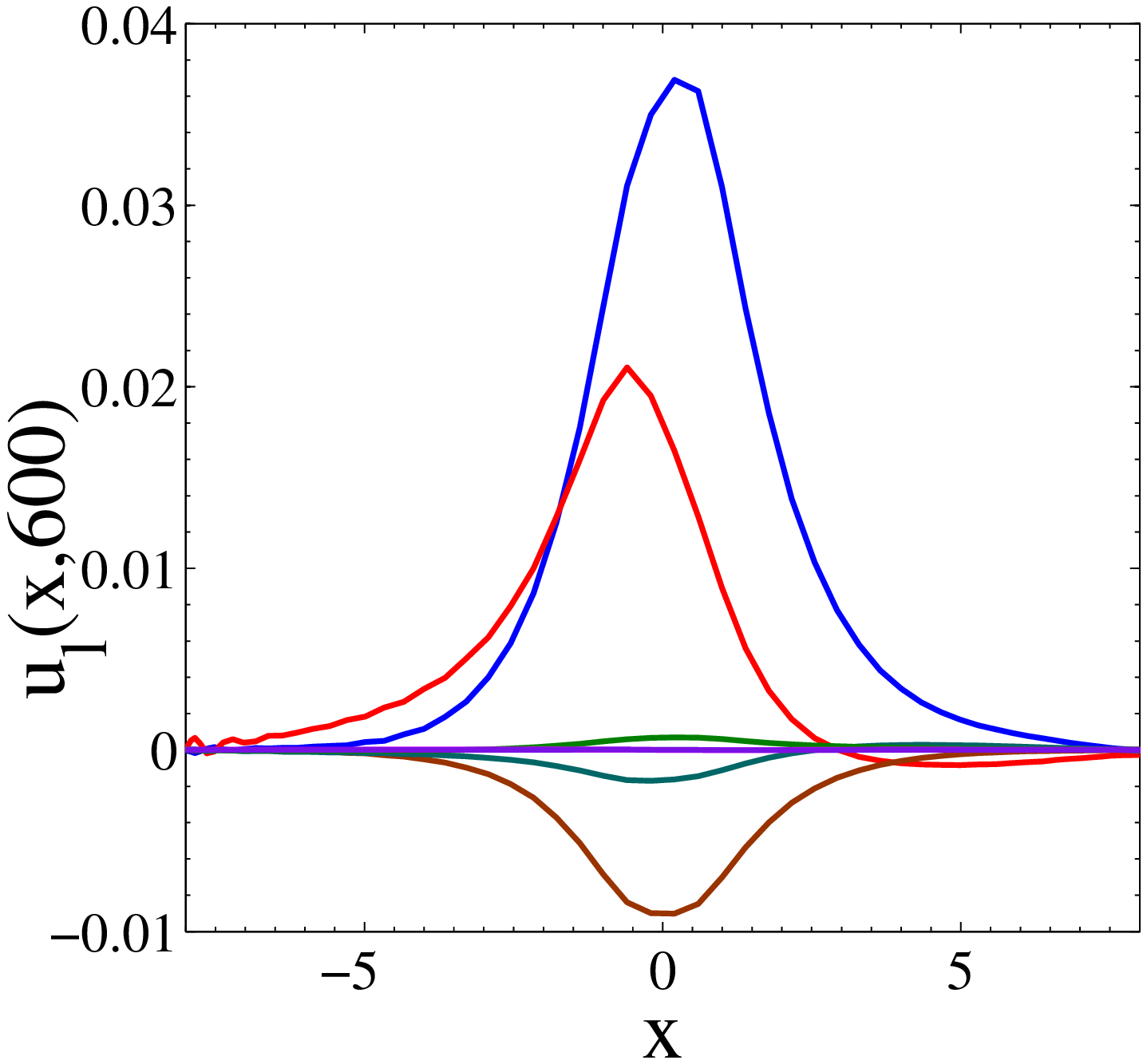}
		\caption{Left: The mean solution  and the standard deviation at $t_{final} = 600$. Right: The first few PC modes for the Hermite chaos at $t_{final} = 600$ and $l = 1,\cdots , 7$.}
	 \label{fig:hermite_sol}
\end{figure}
The boundary conditions are given by
\begin{eqnarray}
\left\lbrace\hat{u}_{l}(L,t) \right\rbrace_{t} &+& \left\lbrace\hat{u}_{l}(L,t) \right\rbrace_{x} = 0, \nonumber \\
\hat{u}_{l}(-L,t) &=& \frac{4}{\sqrt{2\pi}l!}\int_{-\infty}^{\infty}\tan^{-1}\left[ \exp \left\lbrace \frac{-L-x_{0}-\left(\mu+\sigma\gamma\right)t}{\sqrt{1-\left(\mu+\sigma\gamma\right)^2}} \right\rbrace \right]H_{l}(\gamma)\omega(\gamma)d\gamma,  \nonumber
\end{eqnarray}
\noindent
for  $ l = 0, 1, \cdots, N.$

\subsection{Mean and Standard Deviation of the solution}
By definition, the mean solution is given by
\begin{eqnarray}
\overline{u}(x,t,V|V\in \left[\alpha, \beta\right]) &=& \int_{\alpha}^{\beta}f(V;\alpha,\beta,\mu,\sigma)u(x,t,V)dV \nonumber \\
&=& \frac{1}{\sqrt{2\pi}\sigma\left[\Phi\left(\frac{\beta - \mu}{\sigma}\right)-\Phi\left(\frac{\alpha - \mu}{\sigma}\right)\right]}\int_{\alpha}^{\beta}e^{-\frac{\left(V-\mu\right)^2}{2\sigma^{2}}}\sum_{l=0}^{\infty}\hat{u}_{l}(x,t)H_{l}\left(\frac{V-\mu}{\sigma}\right)dV \nonumber \\
&\approx & \hat{u}_{0}(x,t) + \frac{1}{\sqrt{2\pi}\left[\Phi\left(\frac{\beta - \mu}{\sigma}\right)-\Phi\left(\frac{\alpha - \mu}{\sigma}\right)\right]}\sum_{l=1}^{N}\hat{u}_{l}(x,t)\int_{\frac{\alpha-\mu}{\sigma}}^{\frac{\beta-\mu}{\sigma}}e^{-\frac{\gamma^2}{2}}H_{l}(\gamma)d\gamma \nonumber\\
&\approx & \hat{u}_{0}(x,t) + \frac{1}{\sqrt{2\pi}\left[\Phi\left(\frac{\beta - \mu}{\sigma}\right)-\Phi\left(\frac{\alpha - \mu}{\sigma}\right)\right]}\sum_{l=1}^{N}\hat{u}_{l}(x,t)\int_{-\infty}^{\infty}e^{-\frac{\gamma^2}{2}}H_{l}(\gamma)d\gamma \nonumber 
\end{eqnarray}
By using the orthogonality of $H_{l}(\gamma)$
$$\overline{u}(x,t) \approx \hat{u}_{0}(x,t).$$

To calculate standard deviation we first calculate $\overline{u^2(x,t)}$.

\begin{eqnarray}
\overline{u^2(x,t)} &=& \int_{\alpha}^{\beta}f\left(V;\alpha ,\beta ,\mu , \sigma\right)u^{2}(x,t,V)dV \nonumber \\
&=& \int_{\alpha}^{\beta} f\left(V;\alpha ,\beta ,\mu , \sigma \right)\left\lbrace \sum_{l=0}^{\infty}\hat{u}_{l}(x,t)H_{l}\left(\frac{V-\mu}{\sigma}\right) \right\rbrace ^2 dV \nonumber \\
&=& \int_{\alpha}^{\beta} f\left(V;\alpha ,\beta ,\mu , \sigma\right)\hat{u}_{0}^{2}(x,t)H_{0}^{2}\left(\frac{V-\mu}{\sigma}\right) dV + \int_{\alpha}^{\beta} f\left(V;\alpha ,\beta ,\mu , \sigma\right)\hat{u}_{1}^{2}(x,t)H_{1}^{2}\left(\frac{V-\mu}{\sigma}\right)dV\nonumber \\
&+& \cdots   + 2\sum_{m}\sum_{n}\int_{\alpha}^{\beta} f\left(V;\alpha ,\beta ,\mu , \sigma\right)\hat{u}_{m}(x,t)\hat{u}_{n}(x,t)H_{m}\left(\frac{V-\mu}{\sigma}\right)H_{n}\left(\frac{V-\mu}{\sigma}\right)dV \nonumber \\
&\approx & \hat{u}_{0}^{2}(x,t) + \frac{1}{\left[\Phi\left(\frac{\beta - \mu}{\sigma}\right)-\Phi\left(\frac{\alpha - \mu}{\sigma}\right)\right]}\sum_{m=1}^{N}m!\hat{u}_{m}^{2}(x,t),\nonumber 
\end{eqnarray}
and 
$$\overline{u^2(x,t)} \approx \hat{u}_{0}^{2}(x,t) + \frac{1}{\left[\Phi\left(\frac{\beta - \mu}{\sigma}\right)-\Phi\left(\frac{\alpha - \mu}{\sigma}\right)\right]}\left[1!\hat{u}_{1}^{2}(x,t) + 2!\hat{u}_{2}^{2}(x,t) + \cdots + N!\hat{u}_{N}^{2}(x,t)\right].$$
The standard deviation is then given by 
$$\sigma = \sqrt{\overline{u^2(x,t)} - \left\lbrace \overline{u}(x,t)\right\rbrace ^2} \approx \left( \frac{1}{\left[\Phi\left(\frac{\beta - \mu}{\sigma}\right)-\Phi\left(\frac{\alpha - \mu}{\sigma}\right)\right]}\sum_{m = 1}^{N}m!\hat{u}_{m}^{2}(x,t)\right) ^{\frac{1}{2}}.$$

\subsection{Estimation of critical velocity} 
To estimate the critical velocity we use the Eq. (\ref{critical_velo_pdf}). Here the density function $f$ is given by Eq. (\ref{hermite-pdf}). 
By using these two equations we get,
\begin{equation}
\Phi\left(\frac{V_{c}-\mu}{\sigma}\right) = \frac{\bar{u}_{\mathrm{her}}}{2\pi}\Phi\left(\frac{\beta-\mu}{\sigma}\right)+ \left(1-\frac{\bar{u}_{\mathrm{her}}}{2\pi}\right)\Phi\left(\frac{\alpha-\mu}{\sigma}\right).
\label{her_critical}
\end{equation}
\noindent
If $\bar{u}_{\mathrm{her}} = 2\pi$, $\Phi\left(\frac{V_{c}-\mu}{\sigma}\right) = \Phi\left(\frac{\beta-\mu}{\sigma}\right)$ which implies $V_{c} = \beta$, all solutions are \textit{trapped} and if $\bar{u}_{\mathrm{her}} = 0$, $\Phi\left(\frac{V_{c}-\mu}{\sigma}\right) = \Phi\left(\frac{\alpha-\mu}{\sigma}\right)$ which implies $V_{c} = \alpha$, all solutions are \textit{passed}.
For numerical simulations, we consider $\mu = 0.12,\; \sigma = 0.01,\; \alpha = 0.11,\; \beta = 0.13.$ We obtain $\bar{u}_{\mathrm{her}} \approx 3.24341621$ at the final time $t_{f} = 600$. By using Eq. (\ref{her_critical}), we get,
$$\Phi\left(\frac{V_{c}-0.12}{0.01}\right) \approx \frac{3.24341621}{2\pi}\left[\Phi\left(\frac{0.13-0.12}{0.01}\right)-\Phi\left(\frac{0.11-0.12}{0.01} \right)\right] + \Phi\left(\frac{0.11-0.12}{0.01}\right).$$
Based on these results and by using the statistical table \cite{Statistics} we estimate $V_{c} \approx 0.1203$. This estimate is close to what we get from the MC simulations \cite{JungDon2009} and the Legendre chaos in the next section. The obtained value is, however, not accurate; it matches with the value obtained by many MC simulations within  two digits only. Figure \ref{fig:hermite_sol} shows the initial condition, the mean solution, and the standard deviation in the left figure
and the first 7 PC modes ($l = 1, \cdots, 7$) in the right figure. 



\section{Sine-Gordon equation with the Legendre chaos} 

%
%

Suppose that $V \in \left[V_{a},V_{b}\right]=[a, b]$ has the uniform distribution and $$V(\xi) = \frac{b-a}{2}\xi + \frac{a+b}{2}, \qquad \xi \in [-1, 1].$$
We consider  $u(x,t,\xi) = \sum_{l=0}^{N} \hat{u}_{l}(x,t)L_{l}(\xi)$.
By using the orthogonal property of the Legendre polynomials and the initial and boundary conditions we have

\begin{table}[h]
\begin{center}
\begin{tabular}{l} 
$\hat{u}_{l}(x,0) = 2(2l+1)\int_{-1}^{1}\arctan\left\lbrace \exp \left[ \frac{x-x_{0}}{\sqrt{1-V^{2}(\xi)}}\right]  \right\rbrace L_{l}(\xi)d\xi,$\\
$\left[ \hat{u}_{l}(x,0)\right]_{t} = -2(2l+1)\int_{-1}^{1}\frac{V(\xi)}{\sqrt{1-V^{2}(\xi)}}\left\lbrace\frac{\exp\left(\frac{x-x_{0}}{\sqrt{1-V^{2}(\xi)}} \right) }{1+\exp\left(\frac{2\left(x-x_{0}\right)}{\sqrt{1-V^{2}(\xi)}} \right) }\right\rbrace L_{l}(\xi)d\xi ,$\\
$\hat{u}_{l}(-L,t) = 2(2l+1)\int_{-1}^{1}\arctan\left\lbrace \exp \left[ \frac{-L-x_{0}-tV(\xi)}{\sqrt{1-V^{2}(\xi)}}\right]  \right\rbrace L_{l}(\xi)d\xi,$\\
$\left\lbrace \hat{u}_{l}(L,t) \right\rbrace_{t} + \left\lbrace \hat{u}_{l}(L,t)\right\rbrace_{x} = 0,$ \\[-1pc]
\end{tabular}
\end{center}
\end{table}
for $l = 0, 1, \cdots , N.$
Then we get
$$\sum_{l=0}^{N}\left[\left\lbrace\hat{u}_{l}(x,t) \right\rbrace_{tt} - \left\lbrace\hat{u}_{l}(x,t) \right\rbrace_{xx}  \right]L_{l}(\xi) = \left\lbrace \epsilon \delta(x)-1 \right\rbrace  \sin\left(\sum_{l=0}^{N}\hat{u}_{l}(x,t)L_{l}(\xi) \right),   $$
for $l = 0, 1, \cdots , N.$
Using the orthogonal property of the Legendre polynomials, we get 
$$\left\lbrace\hat{u}_{l}(x,t) \right\rbrace_{tt} - \left\lbrace\hat{u}_{l}(x,t) \right\rbrace_{xx} = \left(\frac{2l+1}{2} \right)\left\lbrace \epsilon \delta(x)-1 \right\rbrace \int_{-1}^{1} \sin\left(\sum_{l=0}^{N}\hat{u}_{l}(x,t)L_{l}(\xi) \right)L_{l}(\xi)d(\xi).   $$
\noindent
for $l = 0, 1, \cdots , N.$
Or 
$$\left\lbrace\hat{u}_{l}(x,t) \right\rbrace_{tt} - \left\lbrace\hat{u}_{l}(x,t) \right\rbrace_{xx} = \left(\frac{2l+1}{2} \right)\left\lbrace \epsilon \delta(x)-1 \right\rbrace \phi_{l}(x,t),   $$
\noindent
for $l = 0, 1, \cdots , N $ and $\phi_{l}(x,t) = \int_{-1}^{1} \sin\left(\sum_{l=0}^{N}\hat{u}_{l}(x,t)L_{l}(\xi) \right)L_{l}(\xi)d(\xi).$

\begin{figure}
	\centering
		\includegraphics[width=0.48\textwidth]{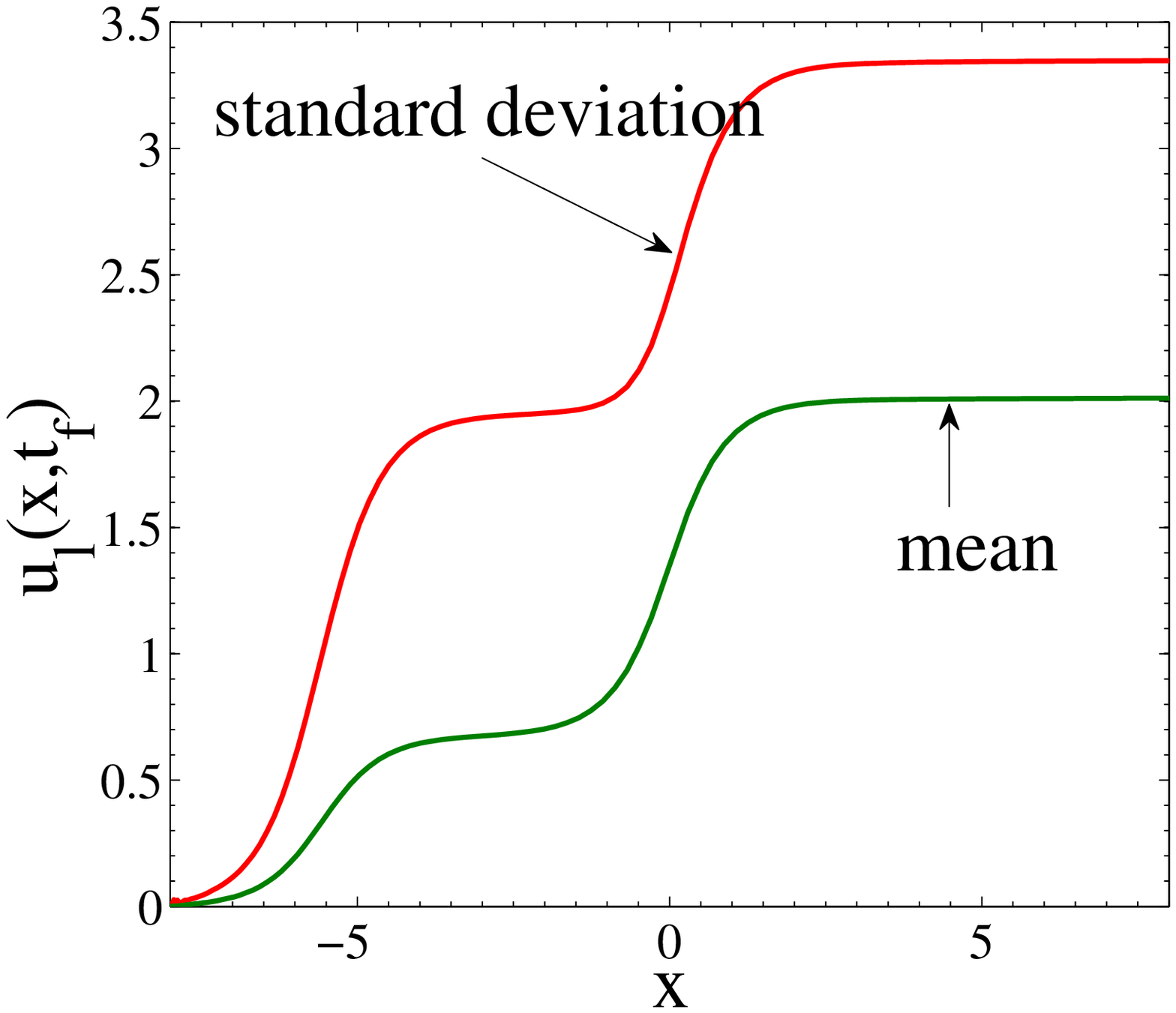}
		\includegraphics[width=0.48\textwidth]{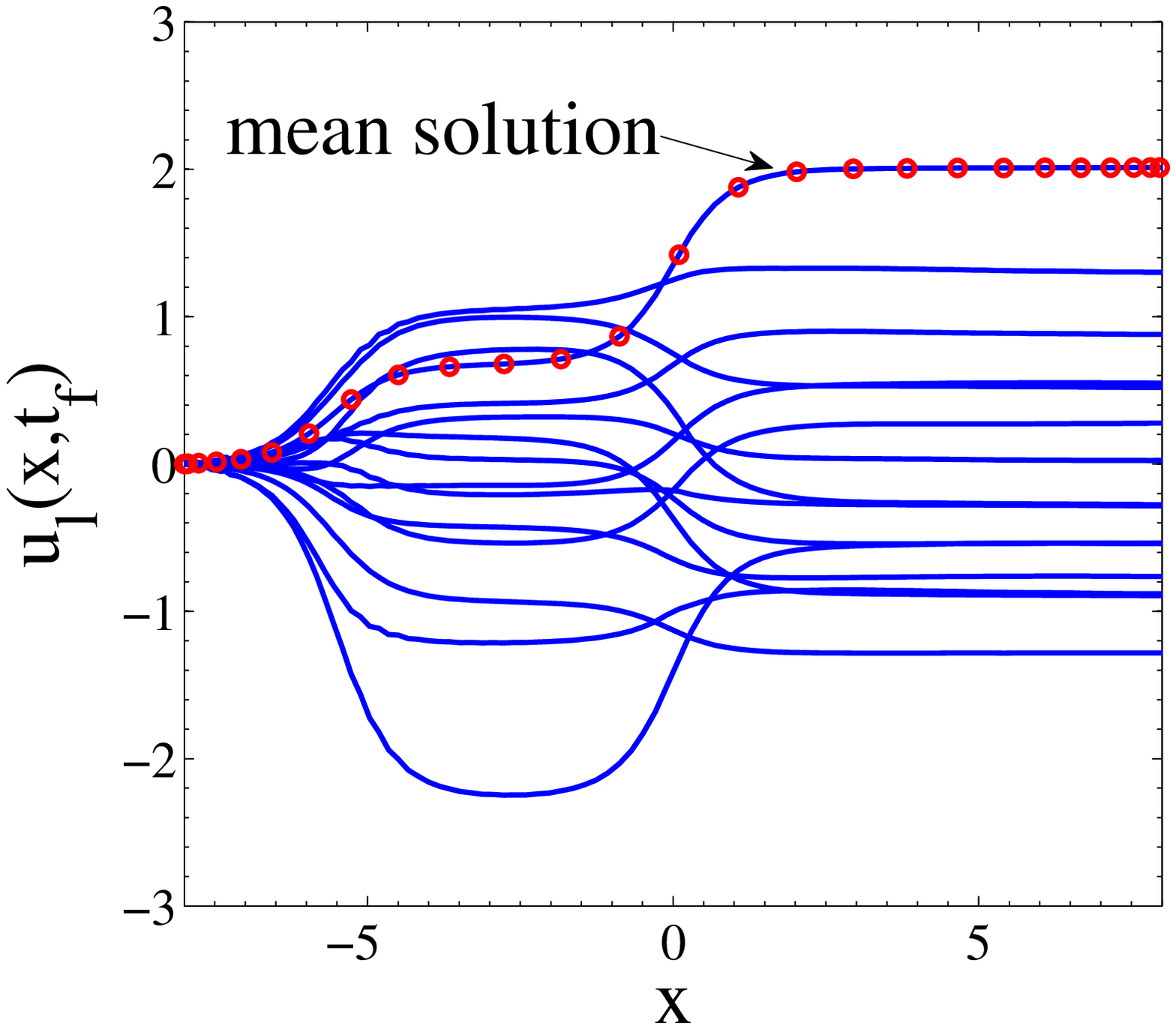}
		\caption{ Left: Mean and Standard deviation of the gPC solution at the final time, $t_{f} = 600$.
		Right: First $15$ gPC modes. Red circles represent the 1st PC mode (mean solution)}
	 \label{pc_critical}
\end{figure}

\begin{figure}
	\centering
		\includegraphics[width=0.6\textwidth]{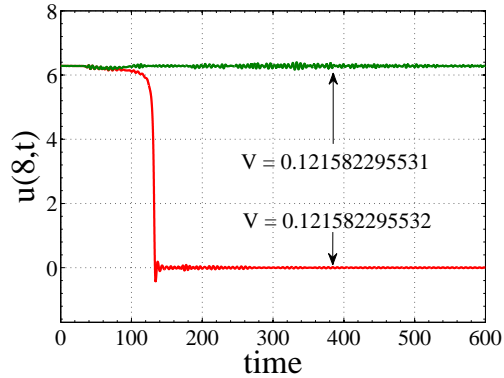}
		\caption{ Solutions at $x = 8$ with time for $V = 0.121582295531$ (sub-critical) and $V = 0.121582295532$ (super-critical). }
	 \label{verifying}
\end{figure}

The discretized equations are given by
\begin{eqnarray}
\hat{u}_{l}^{n+1} &=&  2\hat{u}_{l}^{n} - \hat{u}_{l}^{n-1} + \left(\delta t \right)^2\left[D^2\left(\hat{u}_{l}^{n} \right)  + \left( \frac{2l+1}{2}\right)\left( \epsilon (\mathbf{D}\cdot H)-I\right)\hat{\phi}_{l}^{n}  \right]  , \nonumber \\
\left[ \hat{u}_{l}(x,0)\right]_{t} &=& -2(2l+1)\int_{-1}^{1}\frac{V(\xi)}{\sqrt{1-V^{2}(\xi)}}\left\lbrace\frac{\exp\left(\frac{x-x_{0}}{\sqrt{1-V^{2}(\xi)}} \right) }{1+\exp\left(\frac{2\left(x-x_{0}\right)}{\sqrt{1-V^{2}(\xi)}} \right) }\right\rbrace L_{l}(\xi)d\xi , \nonumber \\
\hat{u}_{l}(-L,t) &=& 2(2l+1)\int_{-1}^{1}\arctan\left\lbrace \exp \left[ \frac{-L-x_{0}-tV(\xi)}{\sqrt{1-V^{2}(\xi)}}\right]  \right\rbrace L_{l}(\xi)d\xi,  \nonumber \\
 \hat{u}_{l}^{n+1}(L,t) &=& \hat{u}_{l}^{n}(L,t) - \left(\frac{\delta t}{\delta x}\right)\left[\hat{u}_{l}^{n}(L,t)-\hat{u}_{l}^{n+1}(L-\delta x,t) \right]   ,\nonumber
\end{eqnarray}
\noindent
for $l = 0, 1, \cdots , N.$
Figure \ref{pc_critical} shows the  PC simulation with $V_{a} = 0.1215$ and $V_{b} = 0.121757$. 
The left figure shows the mean solution and the standard deviation and the right figure the first $15$ PC modes. The 
first mean mode is plotted with the red circle. 
From these results we find the mean solution $\bar{u}_{leg}$, at $x = L$. $\bar{u}_{leg} = 2.011976945534794 $ at $t_{f} = 600$.
By using Eq. (\ref{critical_legendre}), we obtain $V_{c} \approx 0.1215822955316$, which is close to the critical value by the MC simulation, which was 
$0.121582 < V_c <0.121583 $ \cite{JungDon2009}. To verify the obtained result, two direct spectral simulations are conducted with 
$V_1 = 0.121582295531$ and $V_2 = 0.121582295532$.  Figure \ref{verifying} shows that the sub-critical solution is obtained for 
$V_1 = V_c - 5.9999\times 10^{-13}$ and the super-critical solution for $V_2 = V_c + 4.0001 \times 10^{-13}$. This verifies that the obtained 
$V_c$ is as accurate as at least up to $13$ digits. 

\begin{remark}
For the Hermite chaos, the integral for the mean solution is not exact and  we need to rely on the statistical table to get the probability values.
For the Legendre chaos, however, we could find a highly accurate result. The artificial oscillations fluctuating around $0$ (particle-pass) and $2\pi$ (particle-capture) shown in the individual solution at $x = L$ due to the non-exact boundary condition and the Gibbs phenomenon are highly reduced in the mean solution and the analysis shown in Diagram 1 yields an accurate mean solution resulting in finding an accurate critical velocity $V_c$. 
\end{remark}
%
%

\section{sine-Gordon equation with the uncertainty in $\epsilon$ }
Suppose that we have the uncertainty in $\epsilon$, and 
$\epsilon \in [a, b]$ has the uniform distribution. Then
$$\epsilon(\xi) = \frac{b-a}{2}\xi + \frac{a+b}{2}, \qquad \xi \in [-1, 1].$$
\noindent
We consider 
$$
u(x,t,\xi) = \sum_{l=0}^{N} \hat{u}_{l}(x,t)L_{l}(\xi). 
$$  

\begin{table}[h]
\begin{center}
\begin{tabular}{l} 
$\hat{u}_{0}(x,0) = 4\arctan\left\lbrace \exp \left[ \frac{x-x_{0}}{\sqrt{1-V^{2}}}\right]  \right\rbrace ,$\\
$\hat{u}_{m}(x,0) = 0,$ \\
$\left[ \hat{u}_{0}(x,0)\right]_{t} = -\frac{4V}{\sqrt{1-V^{2}}}\left\lbrace\frac{\exp\left(\frac{x-x_{0}}{\sqrt{1-V^{2}}} \right) }{1+\exp\left(\frac{2\left(x-x_{0}\right)}{\sqrt{1-V^{2}}} \right) }\right\rbrace  ,$\\
$\left[ \hat{u}_{m}(x,0)\right]_{t}  = 0,$ \\
$\hat{u}_{0}(-L,t) = 4\arctan\left\lbrace \exp \left[ \frac{-L-x_{0}-tV}{\sqrt{1-V^{2}}}\right]  \right\rbrace,$\\
$\hat{u}_{m}(-L,t) = 0,$ \\
$\left\lbrace \hat{u}_{l}(L,t) \right\rbrace_{t} + \left\lbrace \hat{u}_{l}(L,t)\right\rbrace_{x} = 0,$\\[-1pc]
\end{tabular}
\end{center}
\end{table}

In the similar way as in the previous section we have for $m =  1, \cdots , N.$

$$\sum_{m=0}^{N}\left[\left\lbrace\hat{u}_{m}(x,t) \right\rbrace_{tt} - \left\lbrace\hat{u}_{m}(x,t) \right\rbrace_{xx}  \right]L_{m}(\xi) = \left\lbrace \left( \frac{b-a}{2}\xi + \frac{a+b}{2} \right)\delta( x)-1 \right\rbrace  \sin\left(\sum_{m=0}^{N}\hat{u}_{m}(x,t)L_{m}(\xi) \right).   $$
Using orthogonal property of the Legendre polynomials again, we get 
$$\left\lbrace\hat{u}_{m}(x,t) \right\rbrace_{tt} - \left\lbrace\hat{u}_{m}(x,t) \right\rbrace_{xx} = \left(\frac{2m+1}{2} \right)\left\lbrace\frac{a+b}{2}\delta(x) -1 \right\rbrace \phi_{m}(\xi) + \frac{b-a}{2}\psi_m(\xi),   $$
\noindent
for $m = 0, 1, \cdots , N.$ where 
 \begin{eqnarray} \phi_{m}(x,t) &=& \int_{-1}^{1} \sin\left(\sum_{l=0}^{N}\hat{u}_{l}(x,t)L_{l}(\xi) \right)L_{m}(\xi)d(\xi), \nonumber \\
 \psi_{m}(x,t) &=& \int_{-1}^{1} \xi \sin\left(\sum_{l=0}^{N}\hat{u}_{l}(x,t)L_{l}(\xi) \right)L_{m}(\xi)d(\xi). \nonumber
\end{eqnarray}

The discretized equations are given by,

\begin{table}[h]
\begin{center}
\begin{tabular}{l} 
$\hat{u}_{m}^{n+1} =  2\hat{u}_{m}^{n} - \hat{u}_{m}^{n-1} + \left(\delta t \right)^2\left[D^2\left(\hat{u}_{m}^{n} \right)  + \left( \frac{2m+1}{2}\right)\left\lbrace \left( \frac{a+b}{2} (\mathbf{D}H)-I\right)\phi_{m}^{n} + \frac{b-a}{2}\psi_{m}^{n}\right\rbrace  \right]  ,$\\
\\
$\hat{u}_{0}(-L,t) = 4\arctan\left\lbrace \exp \left[ \frac{-L-x_{0}-tV}{\sqrt{1-V^{2}}}\right]  \right\rbrace ,$\\
\\
$\hat{u}_{m}(-L,t) = 0,$ \\
\\
$ \hat{u}_{m}^{n+1}(L,t) = \hat{u}_{m}^{n}(L,t) - \left(\frac{\delta t}{\delta x}\right)\left[\hat{u}_{m}^{n}(L,t)-\hat{u}_{m}^{n+1}(L-\delta x,t) \right],$\\[-1pc]
\end{tabular}
\end{center}
\end{table}
\noindent
for $m = 0, 1, \cdots , N.$

\subsection{Finding critical amplitude}

%
%

To get the critical value of $\epsilon$, $\epsilon_{c}$, we follow the same procedure as in section 6.2. Since $\epsilon$ is uniformly distributed in $\left[\epsilon_{a}, \epsilon_{b}\right]$, we use Eq. (\ref{critical_legendre}), but in little different way.
Let $\epsilon_{a} \le \epsilon_{c} \le \epsilon_{b}$ and for $\epsilon \in \left[\epsilon_{a}, \epsilon_{c}\right] $, the solutions are \textit{passed} and for $\epsilon \in \left(\epsilon_{c}, \epsilon_{b}\right] $, the solution get \textit{trapped}. So Eq. (\ref{critical_legendre}) takes the form  
\begin{eqnarray}
\bar{u}_{leg} &=& \frac{2\pi}{\epsilon_{b}-\epsilon_{a}}\left[\epsilon_{b}- \epsilon_{c}\right] \nonumber \\
\epsilon_{c} &=& \epsilon_{b} -\left(\epsilon_{b}-\epsilon_{a} \right)\frac{\bar{u}_{leg}}{2\pi}.
\label{critical_eps_leg}
\end{eqnarray}

\begin{figure}
	\centering
		\includegraphics[width=0.6\textwidth]{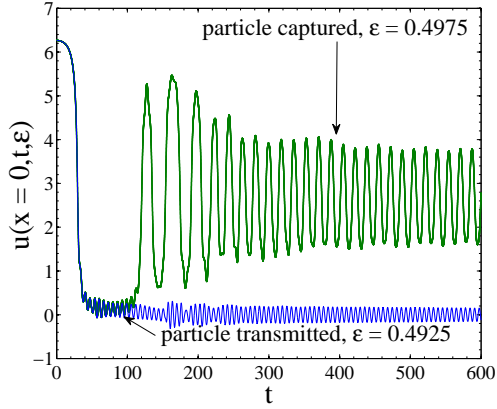}
		\caption{ MC simulations for different $\epsilon$ to get the idea of the position of $\epsilon_{c}$.}
	 \label{sg_epsilon_mc}
\end{figure}

\begin{figure}
	\centering
		\includegraphics[width=0.6\textwidth]{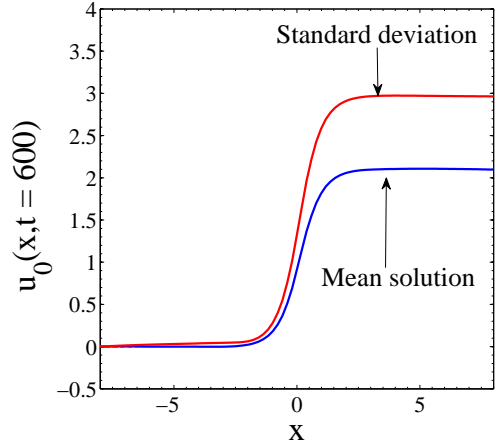}
		\caption{ 
		 Mean and Standard deviation of the PC solution at $T_{f} = 600$ where $\epsilon$ is stochastic variable. }
	 \label{sg_epsilon_critical}
\end{figure}

Figure (\ref{sg_epsilon_mc}) shows some solutions at $x = 0$ with time for  different values of $\epsilon$, $\epsilon = 0.4925, 0.4975$. 
For the captured solution, the solution at $x = 0$ is highly oscillatory but for the transmitted one it is close to zero. 
From the figure we know that $0.4925 < V_c < 0.4975$. 
%
%
\noindent
Suppose we know for $\epsilon_{a} = 0.495$ solution passes through the barrier and for $\epsilon_{b} = 0.4975$ solution is \textit{trapped} by the barrier. By the gPC simulations (Figure \ref{sg_epsilon_critical}) for $t_{f} = 600$ we have $\bar{u}_{leg} = 2.096586$ . Using these in Eq. (\ref{critical_eps_leg})  we get 
$\epsilon_{c}= 0.4958342050637932$.

\section{Conclusion}
In this paper, we proposed an efficient method using the gPC method to find the critical parameter values and the related statistical quantities for the Klein-Gordon and sine-Gordon equations with a point-like singular potential function. By assuming the unknown critical parameter values as a random variable, and expanding the solution in the orthogonal polynomials associated with the random variable, we computed the critical parameter with much less computation than the MC approach. For the Klein-Gordon equation, the critical parameter is associated with the amplitude of the singular potential function. Using the gPC reconstruction of the solution, the invariant point is found to be the critical parameter value. For the sine-Gordon equation, the initial velocity of the soliton solution is used as the random variable. For this case, both the Hermite and Legendre gPC methods are used. Due to the nonlinearity and the $\delta$-function like PDF of the receiving signals, the reconstruction could not be used; the gPC mean was used instead. The Legrendre gPC mean converges quickly, and it efficiently and accurately determines the critical parameter value. 

As our numerical results show, the gPC method offers several advantages in finding the critical value as compared to the MC method. First of all, one can avoid the large number of individual simulations for different parametric values. Our proposed method suggests that, with a small number of MC simulations, it is possible to make the first guess about the critical value. Once the first guess is made, one can apply the gPC method to narrow down or pinpoint the critical value accurately and efficiently. Moreover, regarding the gPC method for the sine-Gordon equation, the stochastic variable appears only in the initial and boundary conditions, and does not appear in the main PDE. To perform the integration related to $V$, one can use the Legendre-Gauss quadrature for the Legendre chaos and Gauss-Hermite quadrature for the Hermite chaos and obtain a fast convergence. In this paper, only $30$ quadrature points are used. Adopting the quadrature rule, one can easily speed up the whole process.

During our study, we also found a limitation of the gPC method regarding the determination of the initial critical velocity or the critical amplitude of the singular potential term for the sine-Gordon equation. That is, we could not use the reconstruction constructed by the gPC modes to find the PDF of the receiving solutions due to the non-linearity and non-separability related to the Gibbs phenomenon of the spectral method. 

Despite this limitation, we find that the gPC method works well for finding the critical values for the singularly perturbed Klein-Gordon and sine-Gordon equations. Since no significant gPC analysis exists for the singularly perturbed PDEs by the Dirac $\delta$-function, our study could be a good resource for this type of research. In our future work, we will apply a similar gPC method for different PDEs, such as the nonlinear Schr\"{o}dinger equations with singular potential terms. We will also consider more general types of uncertainties associated with the singular potential term, for which more than one random variables are involved and the gPC expansion should be in the multi-dimensional random space. To reduce the complexity due to the dimensionality, we will attempt to use the gPC collocation method with the fast sparse grid methods \cite{LiuGaoHesthaven,Dongbin2007, Xiu2005}.  

\vskip .1in
{\bf Acknowledgement:}
\noindent
We are grateful to Bruce Pitman, Gino Biondini, and Emmanuel Lorin for their useful discussions and references.

\vskip .1in
\begin{center}
{\bf Appendices}
\end{center} 

\appendix

{\bf Appendix 1: Hermite polynomials}

%
%
%
%
%
%
%
$H_{n}(x)$ are the Hermite polynomial for $n = 0, 1, 2, 3 \cdots $ and they are orthogonal with respect to the weight function 
$w(x) = e^{-\frac{x^2}{2}}$; 
$$\int_{-\infty}^{\infty}H_{m}(x)H_{n}(x)e^{-\frac{x^2}{2}}dx = \sqrt{2\pi}n!\delta_{mn}.$$
First few Hermite polynomials are 
$H_{0}(x) = 1, H_{1}(x) = x, H_{2}(x) = x^{2} -1,  H_{3}(x) = x^{3} - 3x. $

%
The following few recursion relations are used in the paper: 
\begin{eqnarray}
H_{n+1}(x) &=& xH_{n}(x) - H_{n}^{'}(x). \nonumber \\
H_{n}^{'}(x) &=& nH_{n-1}(x). \nonumber \\
H_{n}(x+y) &=& \sum_{k=0}^{n}C^{n}_{k}x^{k}H_{n-k}(y). \nonumber 
\end{eqnarray}


{\bf Appendix 2: Legendre Polynomials}

Some important properties of the Legendre polynomials are used in the paper.  
Legendre polynomials are defined by the solutions $p_n(x)$  of the following Sturm-Liouville equation: 
\begin{equation} 
\frac{d}{dx}[(1-x^2)\frac{dp_{n}}{dx}]+n(n+1)p_{n}=0, \quad x \in [-1,1], 
\end{equation}
where $p_n(x)$ is the Legendre polynomial of degree $n$.  First few Legendre polynomials are 
 $p_{0}(x)=1,$  $p_{1}(x)=x,$  $p_{2}(x)=\frac{1}{2}(3x^2-1),$ $p_{3}(x)=\frac{1}{2}(5x^3-3x),$ $p_{4}(x)=\frac{1}{8}(35x^4-30x^2+3)$. 
 Also $p_{n}(1)=1$, \ $p_{n}(-1) = 1$ if $n$ is even and $= -1$ if $n$ is odd, also $p_{n}(0)=0$ if $n$ is odd with the orthogonality
 $$
 \int_{-1}^{1}p_{m}(x)p_{n}(x)dx = \frac{2}{2n+1}\delta_{mn}.
 $$
 
The recurrence relation is given by 
\begin{eqnarray}
(n+1)p_{n+1}(x)= (2n+1)xp_{n}(x)-np_{n-1}(x),\\
p_{n+1}^{'}(x) -p_{n-1}^{'}(x)=(2n+1)p_{n}(x).
\label{legendre_recurrence}
\end{eqnarray}

From  Eq. (.2) we have,
$$
(2n+1)xP_{n}(x) = nP_{n-1}(x) + (n+1)P_{n+1}(x). 
$$
Multiplying both sides by $P_{m}$ and integrating from $-1$ to $1$, we get,
$$
\int_{-1}^{1}xP_{n}(x)P_{m}(x)dx = \frac{n}{2n+1}\int_{-1}^{1}P_{n-1}(x)P_{m}(x)dx + \frac{n+1}{2n+1}\int_{-1}^{1}P_{n+1}(x)P_{m}(x)dx.
$$
\noindent
By using the orthogonality of the Legendre polynomials we get,
\begin{equation}
\int_{-1}^{1}xP_{n}(x)P_{m}(x)dx = \left\{\begin{array}{cc} \frac{2(m+1)}{(2m+1)(2m+3)} ,&\mathrm{if} \quad n = m+1 , \\
   \frac{2m}{(2m+1)(2m-1)}, & \mathrm{if} \quad n = m-1,  \\
    0 ,&\mathrm{otherwise}. \end{array} \right.
\end{equation}

\noindent
By using Eq. (.3) we get the followings.\\
For $n$ is odd number: 
\begin{equation}
p^{'}_{n+1}=(2n+1)p_{n}(x) + (2n-3)p_{n-2}(x) + \cdots + 7p_{3}(x) + 3p_{1}(x).
\label{odd_derivatives}
\end{equation}
For $n$ is even number: 
\begin{equation}
p^{'}_{n+1}=(2n+1)p_{n}(x) + (2n-3)p_{n-2}(x) + \cdots + 9p_{4}(x) + 5p_{2}(x) +p_{0}(x).
\label{even_derivatives}
\end{equation}

\noindent
The derivative of Legendre polynomials as the finite linear sum of Legendre polynomials are used in the paper and 
some examples are as follows:
\noindent
For $n = 6$, 
$$p_{7}^{'}(x) = 13p_{6}(x)+9p_{4}(x)+5p_{2}(x) +p_{0}(x).$$ 
And for $n = 7$,  
$$p_{8}^{'}(x) = 15p_{7}(x)+11p_{5}(x)+7p_{3}(x) +3p_{1}(x).$$
By orthogonal property $\int_{-1}^{1}p^{'}_{7}(x)p^{'}_{8}(x)dx = 0$. So in general if $m$ is even and $n$ is odd or $n$ is even and $m$ is odd,
\begin{equation}\int_{-1}^{1}p^{'}_{m}(x)p^{'}_{n}(x)dx = 0.\end{equation} \\
If both $m$ and $n$ are even numbers and $m \ge n$ then by using the orthogonal property of Legendre polynomials and from Eq. (\ref{even_derivatives}) we have,
\begin{eqnarray}
\int_{-1}^{1}p^{'}_{m}(x)p^{'}_{n}(x)dx &=& (2n-1)^2\frac{2}{2(n-1)+1} + (2n-5)^2\frac{2}{2(n-3)+1} + \nonumber  \\
&&\cdots + 7^2\frac{2}{2\cdot 3+1}+3^2\frac{2}{2\cdot 1+1}\nonumber \\
&=&2[3+7+\cdots +(2n-1)]=n(n+1).\nonumber 
\end{eqnarray} 
If both $m$ and $n$ are odd numbers and $m \ge n$ then the orthogonal property and Eq. (\ref{odd_derivatives}) yield, \begin{equation}\int_{-1}^{1}p^{'}_{m}(x)p^{'}_{n}(x)dx = 2[1+5+9+\cdots + (2n-1)]=n(n+1).\end{equation} \\ \\

\end{document}